\documentclass[a4paper,11pt]{article}

\usepackage{amsfonts,amssymb,amsbsy,amsmath,amsthm,mathrsfs,enumerate,verbatim}

  \usepackage{graphics} %% add this and next lines if pictures should be in esp format
  \usepackage{epsfig} %For pictures: screened artwork should be set up with an 85 or 100 line screen
\usepackage{graphicx}  \usepackage{epstopdf}%This is to transfer .eps figure to .pdf figure; please compile your paper using PDFLeTex or PDFTeXify.
 \usepackage[colorlinks=true]{hyperref}

\hypersetup{urlcolor=blue, citecolor=red}

\usepackage{amsfonts}
\usepackage{amsmath}
\usepackage{amssymb}
\usepackage{graphicx}
\usepackage{epsf}
\usepackage{epsfig}

\usepackage{color}
\usepackage{afterpage}
\usepackage{verbatim}
\usepackage{times}
\newtheorem{theo}{Theorem}

\newtheorem{theorem}{Theorem}[section]

\newtheorem{lemma}[theorem]{Lemma}

\theoremstyle{definition}
\newtheorem{remark}[theorem]{Remark}

\newcommand{\bel}{\begin{equation} \label}
\newcommand{\ee}{\end{equation}}

\def\beq{\begin{equation}}
\def\eeq{\end{equation}}
\newcommand{\bea}{\begin{eqnarray}}
\newcommand{\eea}{\end{eqnarray}}
\newcommand{\beas}{\begin{eqnarray*}}
\newcommand{\eeas}{\end{eqnarray*}}

{

 \usepackage{graphicx,color}
 \definecolor{mygreen}{cmyk}{1,0,1,0.1}

\usepackage{epsf}
\usepackage{epstopdf}

\usepackage{pdfsync}

\graphicspath{% directories to search the graphics files
  {Figures/}
}

\begin{document}

%\title{Optimization approach for determining the conductivity in a
%  waveguide}
\title{An explicit $P_1$ finite element scheme for
  Maxwell's equations with constant permittivity in a boundary neighborhood}

\author{L. Beilina  \thanks{
Department of Mathematical Sciences, Chalmers University of Technology and
 University of Gothenburg, SE-42196 Gothenburg, Sweden, e-mail: \texttt{\
larisa@chalmers.se}}
\and
V. Ruas \thanks{Institut Jean Le Rond d'Alembert, UMR 7190 CNRS - Sorbonne Universit\'e, F-75005, Paris, France, e-mail: \texttt{\
vitoriano.ruas@upmc.fr}}
}

\date{}

\maketitle

%The abstract of your paper
\begin{abstract}

  This paper is devoted to the complete convergence study of the approximation of Maxwell's equations in terms of the sole electric field, 
by means of standard linear finite elements for the space discretization, combined with a well-known explicit finite-difference scheme for the time discretization. The analysis applies to the particular case where the electric permittivity has a constant value outside a sub-domain, whose closure does not intersect the boundary of the domain where the problem is defined. Optimal convergence results are derived in natural norms under reasonable assumptions, provided a classical CFL condition  holds.     

\end{abstract}

\maketitle

%%%%%%%%%%%%%%%%%%%%%%%%%%%%%%%%%%%%%%%%%%%%%%%%%%%%%%%%%%%%%%%%%%%%%%%%%%%%%%%%%%%%%%%%%%%%%%%%%%%%%%%%%%%%%%%%%%%%%%%%%%%%%%%%%%%%%%%%%%%%%%%%%%%%%%%%%%%%%%%%%%%%%%%%%%%%%%%%%%%
%%%%%%%%%%%%%%%%%%%%%%%%%%%%%%%%%%%%%%%%%%%%%%%%%%%%%%%%%%%%%%%%%%%%%%%%%%%%%%%%%%%%%%%%%%%%%%%%%%%%%%%%%%%%%%%%%%%%%%%%%%%%%%%%%%%%%%%%%%%%%%%%%%%%%%%%%%%%%%%%%%%%%%%%%%%%%%%%%%%

\section{Motivation}

The purpose of this article is to present a convergence analysis of an explicit $P1$ finite element solution 
scheme of hyperbolic Maxwell's equations for the electric field with
constant dielectric permittivity in a neighborhood of the boundary of the computational domain. The technique of analysis is inspired by those developed in 
\cite{Ruas, Thomee}.\\

The standard continuous $P1$ FEM is a tempting possibility to solve Maxwell's equations, owing to its simplicity. 
It is well known however that, for different reasons, this method is not always well suited for this purpose. The first reason is that in general the natural function space for the electric field is not the Sobolev space ${\bf H}^1$, but rather in the space ${\bf H}(curl)$. Another issue difficult to overcome with continuous Lagrange finite elements is the prescription of the zero tangential-component boundary conditions for the electric field, which hold in many important applications. All this  motivated the proposal by N\'ed\'elec about four decades ago of a family of ${\bf H}(curl)$-conforming methods to solve these equations (cf. \cite{Nedelec}). These methods are still widely in use, as much as other approaches well adapted to such specific conditions (see e.g. \cite{Assous}, \cite{CiarletJr} and \cite{AMIS}). A comprehensive description of finite element methods for Maxwell's equations can be found in \cite{Monk}. \\
\indent There are situations however in which the $P1$ finite element method does provide an inexpensive and reliable way to solve the Maxwell's equations. In this work we consider one of such cases, characterized by the fact that the electric permittivity is constant in a neighborhood of the whole boundary of the domain of interest. This is because, at least in theory, whenever the electric permittivity is constant, the Maxwell's equations simplify into as many wave equations as the space dimension under consideration. More precisely we show here that, in such a particular case, a space discretization by means of conforming linear elements, combined with a straightforward explicit finite-difference scheme for the time discretization, gives rise to optimal approximations of the electric field, as long as a classical CFL condition is satisfied. \\    
\indent Actually this work can be viewed as both a continuation and the completion of studies presented in
\cite{BG, cejm} for a combination a the finite difference method in a sub-domain with constant permittivity with the finite element method in the complementary sub-domain. As pointed out above, the Maxwell's equations reduces to the wave equation in the former case. Since the analysis of finite-difference methods for this type of  equation is well established, only an explicit $P1$ finite element scheme for Maxwell's equations is analyzed in this paper. \\
\indent In \cite{BG, cejm} a stabilized domain-decomposition finite-element/finite-difference approach for the
solution of the time-dependent Maxwell's system for the electric field was proposed and numerically verified. 
In these works \cite{BG, cejm} different manners to handle a divergence-free
condition in the finite-element scheme were considered. The main idea behind 
the domain decomposition methods in \cite{BG, cejm} is that a rectangular  
computational domain is decomposed into two sub-domains, in which 
two different types of discretizations are employed, namely, the finite-element domain 
in which a classical $P1$ finite element discretization is used, and the finite-difference domain, 
in which the standard five- or seven-point finite difference scheme 
is applied, according to the space dimension. The finite element domain 
lies strictly inside the finite difference domain, in such a way that both domains overlap in two layers
of structured nodes. First order absorbing boundary conditions
\cite{EM} are enforced on the boundary of the computational domain, i.e. on the outer boundary of the finite-difference
domain. In \cite{BG, cejm} it was assumed that the dielectric
permittivity function is strictly positive and has a constant value in
the overlapping nodes as well as in a neighborhood of the boundary of
the domain. An explicit scheme was used both in the finite-element
and finite-difference domains. \\%so that P1 finite element method was
%applied in the finite element domain.
\indent We recall that for a stable finite-element
solution of Maxwell's equation divergence-free edge elements are the
most satisfactory from a theoretical point of view ~\cite{Nedelec,
  Monk}.  However, the edge elements are less attractive for solving 
time-dependent problems, since a linear system of equations should
be solved at every time iteration.  In contrast, $P1$ elements can be
efficiently used in a fully explicit finite element scheme with lumped
mass matrix \cite{delta, joly}. On the other hand it is also well known that the 
numerical solution of Maxwell's equations with nodal finite elements can 
result in unstable spurious solutions \cite{MP, PL}. Nevertheless a number of techniques 
are available to remove them, and in this respect we refer for
example to \cite{Jiang1, Jiang2, Jin, div_cor, PL}. 
In the current work, similarly to \cite{BG, cejm}, the spurious
solutions are removed  from the finite element scheme by adding the
divergence-free condition to the model equation for the electric
field.  Numerical tests given in \cite{cejm} demonstrate that 
spurious solutions are removable indeed, in case an explicit $P1$ finite-element
solution scheme is employed.\\
%The first stable time-domain decomposition method for solution of
%Maxwell's equations was proposed In ~\cite{GEMS1, RB1, RB2}.  This
%method combined finite difference time-domain method (FDTD) of
%\cite{Yee} on the structured part of the mesh with tetrahedral edge
%elements on the unstructured part. In works \cite{GEMS1, RB1, RB2} a
%finite element method was implemented using edge elements of
%N\'ed\'elec \cite{Nedelec} on a hexahedral mesh for
%$H(\mbox{curl})$-conforming discretization of the electric field. In
%~\cite{GEMS1,RB1,RB2} implicit time-stepping was used inside finite
%element domain to obtain stability of the whole hybrid scheme in time.
\indent Efficient usage of an explicit $P1$ finite-element scheme for the 
solution of coefficient inverse problems (CIPs), in the particular context described above 
was made evident in \cite{BK}.
 In many algorithms aimed at solving electromagnetic CIPs, a qualitative collection of
experimental measurements is necessary on the boundary of a 
computational domain, in order to determine the dielectric permittivity function
therein. In this case, in principle the numerical solution of the time-dependent
Maxwell's equations is required in the entire space
$\mathbb{R}^{3}$ (see e.g. \cite{BK, BCN, BTKM, MalmbergBeilina1,
  MalmbergBeilina2, Malmberg}, but instead it can be more    
efficient to consider Maxwell's equations with a constant dielectric
permittivity in a neighborhood of the boundary of a 
computational domain.  
The explicit $P1$ finite-element scheme considered in this work was numerically
tested in the solution of the time-dependent Maxwell's system in both two- and three-dimensional 
geometry (cf. \cite{cejm}). It was also combined with a few algorithms to solve different CIPs for
determining the dielectric permittivity function in connection with the time-dependent
Maxwell's equations, using both simulated and experimentally generated data (see \cite{BCN, BTKM, MalmbergBeilina1, MalmbergBeilina2, Malmberg}).
In short, the formal reliability analysis of such a method conducted in this work, corroborates the previously observed adequacy of this numerical approach.\\

An outline of this paper is as follows: In Section 2 we describe in
detail the model problem being solved, and give its equivalent
variational form. In Section 3 we set up the discretizations of the
model problem in both space and time. Section 4 is devoted to the
stability analysis of the explicit scheme considered in the previous
section, and Section 5 to the corresponding consistency study. Next we
combine the results of the two previous sections to prove error
estimates in Section 6. Underlying convergence results under the very
realistic assumption that the time step varies linearly with the mesh
size as the meshes are refined are thus established. In Section 7 
we present a numerical validation of our scheme. Finally we
conclude in Section 8 with a few comments on the whole work.

\section{The model problem}

The Maxwell's equations for the electric field ${\bf e}=(e_1, e_2, e_3)$ in a bounded domain $\Omega$ of $\Re^3$ with boundary $\partial \Omega$ that we study in this work is 
as follows. First we consider that $\Omega = \bar{\Omega}_{in} \cup \Omega_{out}$, where $\Omega_{in}$ is an interior open set whose boundary does not intersect 
$\partial \Omega$ and $\Omega_{out}$ is the complementary set of $\bar{\Omega}_{in}$ with respect to $\Omega$. ${\bf n}$ being the unit outer normal vector on $\partial \Omega$ we denote by $\partial_n(\cdot)$ the outer normal derivative of a field on $\partial \Omega$. Now in case ${\bf e}$ satisfies absorbing boundary conditions, given ${\bf e}_0 \in [H^1(\Omega)]^3$ and ${\bf e}_1 \in {\bf H}(div,\Omega)$ satisfying $\nabla \cdot (\varepsilon {\bf e}_0) = \nabla \cdot (\varepsilon {\bf e}_1) = 0$ where $\varepsilon$ is the electric permittivity. $\varepsilon$ is assumed to belong to $C^{2,\infty}(\bar{\Omega})$ and to fulfill $\varepsilon \equiv 1 $ in $\Omega_{out}$ and $\varepsilon \geq 1$. Incidentally, throughout this article we denote the standard semi-norm of $C^m(\bar{\Omega})$ by $| \cdot |_{m,\infty}$ for $m >0$ and the standard norm of $C^{0}(\bar{\Omega})$ by $\| \cdot \|_{0,\infty}$. \\
\indent In doing so, the problem to solve is:
\begin{equation}\label{eq1}
  \begin{array}{ll}
    \varepsilon \partial_{tt} {\bf e} + \nabla \times  \nabla \times {\bf e}  = {\bf 0} & \mbox{ in } \Omega \times (0, T), \\
    {\bf e}(\cdot,0) = {\bf e}_0(\cdot), \mbox{ and } \partial_t{\bf e}(\cdot,0) = {\bf e}_1(\cdot) & \mbox{ in } \Omega, \\
    \partial _{n} {\bf e} = - \partial_t {\bf e} & \mbox{ on } \partial \Omega \times (0,T), \\
    \nabla \cdot (\varepsilon {\bf e}) = {\bf 0} & \mbox{ in } \Omega.
  \end{array}
\end{equation}

\begin{remark} The study that follows also apply to the case where boundary conditions other than absorbing boundary conditions $\partial_n {\bf e}=-\partial_t {\bf e}$ are prescribed, for which the same qualitative results hold. As pointed out in Section 1, the choice of the latter here was motivated
by the fact that they correspond to practical situations addressed in \cite{BCN, BTKM, MalmbergBeilina1,MalmbergBeilina2, Malmberg}. %The same situation applies to the assumptions on $\varepsilon$. 
\rule{2mm}{2mm}
\end{remark} 
\begin{remark} The assumption that $\varepsilon$ attains a minimum in an outer layer is not essential for our numerical method to work. However, as far as we can see, it is a condition that guarantees optimal convergence results. In the final section a more elaborated discussion on this issue can be found. \rule{2mm}{2mm}
\end{remark} 

\subsection{Variational form} 
Let us denote the standard inner product of $[L^2(\Omega)]^M$ by $(\cdot,\cdot)$ for $M \in \{1,2,3\}$ and the corresponding norm by $\parallel \{\cdot\} \parallel$. Similarly we denote by $(\{\cdot\},\{\cdot\})_{\partial \Omega}$ the standard inner product of $[L^2(\partial \Omega)]^M$ and the associated norm by $\| \{\cdot\} \|_{\partial \Omega}$.  Further, for a given negative function $\omega \in L^{\infty}(\Omega)$ we introduce the weighted $L^2(\Omega)$-semi-norm 
$\| \{\cdot\} \|_{\omega}:=\sqrt{\int_{\Omega} |\omega| |\{\cdot\}|^2 d{\bf x}}$, which is actually a norm if $\omega \neq 0$ everywhere in $\bar{\Omega}$. We also 
introduce, the notation $(A,B)_{\omega}:= \int_{\Omega} \omega A \cdot B d{\bf x}$ for two fields $A,B$ which are square integrable in $\Omega$. Notice that if  
$\omega$ is strictly positive this expression defines an inner product associated with the norm $\| \{\cdot\} \|_{\omega}$.\\
Then requiring that ${\bf e}_{|t=0} = {\bf e}_0$ and $\{\partial_t{\bf e}\}_{|t=0} = {\bf e}_1$, we write for all $ {\bf v} \in [H^1(\Omega)]^3$,
\begin{equation}\label{eq2}
  %\begin{split}
 \left (\partial_{tt} {\bf e},{\bf v} \right )_{\varepsilon} + (\nabla {\bf e},\nabla {\bf v})  + (\nabla \cdot \varepsilon {\bf e}, \nabla \cdot {\bf v}) - 
(\nabla \cdot {\bf e}, \nabla \cdot {\bf v}) + (\partial_t{\bf e}, {\bf v})_{\partial \Omega}  = 0 \;\forall t \in (0, T).
 % \end{split}
\end{equation}

Problem (\ref{eq2}) is equivalent to Maxwell's equations  \eqref{eq1}. Indeed 
integrating by parts  \eqref{eq2}, for all $ {\bf v} \in [H^1(\Omega)]^3$ we get,
\begin{equation}\label{eq3}
\begin{array}{l}
\left ( \partial_{tt} {\bf e},{\bf v} \right )_{\varepsilon} + ( \nabla \times \nabla \times {\bf e}, {\bf v}) - (\nabla \nabla \cdot \varepsilon {\bf e},{\bf v}) \\
+ \left(\partial_n{\bf e} + \partial_t {\bf e},{\bf v} \right)_{\partial \Omega} + (\nabla \cdot \varepsilon {\bf e} 
- \nabla \cdot {\bf e}, {\bf v} \cdot {\bf n} )_{\partial \Omega}  = 0
\end{array}
\end{equation}
Noting that $\varepsilon=1$ on $\partial \Omega$  we get
\begin{equation}\label{eq3bis}
  \begin{array}{ll}
    \varepsilon \partial_{tt} {\bf e}
    +  \nabla \times \nabla \times {\bf e} - \nabla (\nabla \cdot \varepsilon {\bf e}) = 0 & \mbox{ in } \Omega \times (0, T),\\
    \partial _{n} {\bf e} = - \partial_t {\bf e} & \mbox{ on } \partial \Omega \times (0,T). 
  \end{array}
\end{equation}
This implies that $\nabla \cdot (\varepsilon {\bf e})=0$. Indeed, let $\tilde{{\bf e}}$ be the unique solution of the Maxwell's equations
\begin{equation}\label{eq4}
  \begin{array}{ll}
    \varepsilon \partial_{tt} \tilde{\bf e}  +
    \nabla \times  \nabla \times \tilde{\bf e}   = 0 & \mbox{ in }\Omega \times (0, T), \\
    \tilde{\bf e}(\cdot,0) = {\bf e}_0(\cdot) \mbox{ and } \partial_t \tilde{\bf e}(\cdot,0) = {\bf e}_1(\cdot) & \mbox{ in } \Omega, \\
    \partial_{n} \tilde{\bf e} = - \partial_t \tilde{\bf e} & \mbox{ on } \partial \Omega \times (0,T), \\
     \nabla \cdot ( \varepsilon \tilde{\bf e}) = 0 & \mbox{ in } \Omega.
  \end{array}
\end{equation}

Using the well-known operator identity $-\nabla^2 (\cdot)= \nabla \times \nabla \times (\cdot) - \nabla \nabla \cdot (\cdot)$, 
$\bar{{\bf e}} := {\bf e} - \tilde{{\bf e}}$ is easily seen to fulfill :
\begin{equation}\label{eq5}
  \begin{array}{ll}
    \varepsilon \partial_{tt} \bar{\bf e}  - \nabla^2 \bar{{\bf e}} 
    - \nabla \nabla \cdot (\varepsilon - 1)  \bar{{\bf e}}   = {\bf 0} & \mbox{ in } \Omega \times (0, T), \\
    \bar{{\bf e}}(\cdot,0) = {\bf 0} \mbox{ and } \partial_t \bar{\bf e}(\cdot,0) = {\bf 0} & \mbox{ in } \Omega, \\
   \partial_{n} \bar{{\bf e}} = - \partial_t \bar{\bf e} & \mbox{ on } \partial \Omega \times (0,T).
  \end{array}
\end{equation}
Now we multiply both sides of (\ref{eq5}) by ${\bf v} \in [H^1(\Omega)]^3$ and integrate the resulting relation 
in $\Omega$. Since $(\nabla \cdot (\varepsilon-1) \bar{\bf e}, {\bf v} \cdot {\bf n} )_{\partial \Omega}=0$, after integration by parts we obtain:
\begin{equation}\label{aux0}
\begin{array}{ll}
  (\partial_{tt} \bar{\bf e}, {\bf v})_{\varepsilon} +
   (\nabla  \bar{{\bf e}},\nabla {\bf v})  + (\nabla \cdot \bar{{\bf e}},\nabla \cdot {\bf v})_{\varepsilon - 1} + (\partial_t \bar{\bf e},{\bf v})_{\partial \Omega} = 0 & 
	\forall t \in (0,T] \\
	& \\
    \bar{{\bf e}}(\cdot,0) = {\bf 0} \mbox{ and } \partial_t \bar{\bf e}(\cdot,0) = {\bf 0} & \mbox{ in } \Omega. 
  \end{array}
\end{equation}
Next we take ${\bf v} = \partial_t \bar{\bf e}$ in (\ref{aux0}) and integrate the resulting relation in $(0,t)$ for $t \in (0,T]$. Using the zero initial conditions satisfied by $\bar{{\bf e}}$ and $\partial_t \bar{\bf e}$, we easily obtain:
\begin{equation}\label{aux1}
  \parallel \partial_t \bar{\bf e} \parallel^2_{\varepsilon} +
   \parallel \nabla  \bar{{\bf e}} \parallel^2  +  \parallel \nabla \cdot \bar{{\bf e}} \parallel^2_{\varepsilon - 1} 
	+ \int_0^t \parallel \partial_t \bar{\bf e}(\cdot,s) \parallel_{\partial \Omega}^2 \; ds = 0. 
\end{equation}
We readily infer from (\ref{aux1}) that $\bar{{\bf e}} \equiv 0$ and hence ${\bf e}$ is the solution of Maxwell's equations \eqref{eq1}.\\

\section{Space-time discretization}

Henceforth, for the sake of simplicity, we assume that $\Omega$ is a polyhedron. 

\subsection{Space semi-discretization} 
Let $V_h$ be the usual $P1$ FE-space of continuous functions related to a mesh ${\mathcal T}_h$ fitting $\Omega$, consisting of tetrahedrons with maximum edge length $h$, belonging to a quasi-uniform family of meshes (cf. \cite{Ciarlet}). Each element $K \in {\mathcal T}_h$ is to be understood as a closed set. \\
Setting ${\bf V}_h := [V_h]^3$ we define ${\bf e}_{0h}$ (resp. ${\bf e}_{1h}$) to be the usual ${\bf V}_h$-interpolate of 
${\bf e}_0$ (resp. ${\bf e}_1$). Then the space semi-discretized problem to solve is \\

\emph{Find }${\bf e}_{h}\in {\bf V}_{h}$ \emph{ such  that $\forall {\bf v} \in {\bf V}_{h}$ }
\begin{equation}\label{eq6}
  \begin{array}{l}
    \left ( \partial_{tt} {\bf e}_{h}, {\bf v} \right )_{\varepsilon} + (\nabla {\bf e}_h,\nabla {\bf v})+ (\nabla \cdot[\varepsilon {\bf e}_h], \nabla \cdot {\bf v} ) 
    -(\nabla \cdot {\bf e}_h, \nabla \cdot {\bf v}) + (\partial_t {\bf e}_{h}, {\bf v})_{\partial \Omega} = 0,  \\
		\\
   {\bf e}_h(\cdot,0) = {\bf e}_{0h}(\cdot) \mbox{ and } \partial_t {\bf e}_h(\cdot,0) = {\bf e}_{1h}(\cdot) \mbox{ in } \Omega.  
  \end{array}
\end{equation}

\subsection{Full discretization}

To begin with we consider a natural centered time-discretization scheme to solve \eqref{eq6}, namely: Given a number $N$ of time steps we define the time increment $\tau := T/N$. Then we approximate ${\bf e}_h(k\tau)$
by ${\bf e}_h^k$ for $k=1,2,...,N$ according to the following FE scheme for $k=1,2,\ldots,N-1$:
\begin{equation}\label{eq7consist}
  \begin{array}{l}
    \displaystyle \left(\frac{{\bf e}_h^{k+1} - 2 {\bf e}_h^k + {\bf e}_h^{k-1}}{\tau^2}, {\bf v} \right)_{\varepsilon} + (\nabla {\bf e}_h^k, \nabla {\bf v}) 
     + (\nabla \cdot \varepsilon {\bf e}_h^k, \nabla \cdot {\bf v}) - (\nabla \cdot {\bf e}_h^k, \nabla \cdot {\bf v}) \\
		 + \displaystyle \left(\frac{{\bf e}_h^{k+1} - {\bf e}_h^{k-1}}{2\tau}, {\bf v} \right)_{\partial \Omega} = 0 \; \forall {\bf v} \in {\bf V}_h,\\
		\\
     {\bf e}_h^0 = {\bf e}_{0h} \mbox{ and } {\bf e}_h^1 = {\bf e}_h^0 + \tau {\bf e}_{1h} \mbox{ in } \Omega.
  \end{array}
\end{equation}

Owing to its coupling with ${\bf e}_h^{k}$ and ${\bf e}_h^{k-1}$ on the left hand side of \eqref{eq7consist}, ${\bf e}_h^{k+1}$ cannot be determined 
explicitly by \eqref{eq7consist} at every time step. In order to enable an explicit solution we resort to the classical mass-lumping technique. We recall that for  
a constant $\varepsilon$ this 
consists of replacing on the left hand side the inner product $( {\bf u},{\bf v})_{\varepsilon}$ (resp. $({\bf u},{\bf v})_{\partial \Omega}$) by an inner product 
$({\bf u},{\bf v})_{\varepsilon.h}$ (resp. $({\bf u},{\bf v})_{\partial \Omega,h}$), using the trapezoidal rule to compute the integral of $\int_K \varepsilon {\bf u}_{|K} \cdot {\bf v}_{|K}d{\bf x}$ (resp. $\int_{K \cap \partial \Omega} {\bf u}_{|K} \cdot {\bf v}_{|K}dS$), for every element $K$ in ${\mathcal T}_h$, where ${\bf u}$ stands for 
${\bf e}_h^{k+1} - 2 {\bf e}_h^k + {\bf e}_h^{k-1}$ (resp. ${\bf e}_h^{k+1} - {\bf e}_h^{k-1}$). It is well-known that in this case the matrix associated with 
$(\varepsilon {\bf e}_h^{k+1},{\bf v})_h$ (resp. $({\bf e}_h^{k+1},{\bf v})_{\partial \Omega,h}$) for ${\bf v} \in {\bf V}_h$, is a diagonal matrix. In our case $\varepsilon$ is not constant, but the same property will hold if we replace in each element $K$ the integral of 
$\varepsilon {\bf u}_{|K} \cdot {\bf v}_{|K}$ in a tetrahedron $K \in {\mathcal T}_h$ or of ${\bf u}_{|F} \cdot {\bf v}_{|F}$ in a face 
$F \subset \partial \Omega$ of a certain tetrahedron $K \in {\mathcal T}_h$ as follows:

\[ 
\begin{array}{l}
\int_K \varepsilon {\bf u}_{|K} \cdot {\bf v}_{|K} d{\bf x} \approx \varepsilon(G_K)  \displaystyle volume(K) \sum_{i=1}^4 \frac{{\bf u}(S_{K,i}) \cdot {\bf v}(S_{K,i})}{4}  \\
\int_{F \subset \partial \Omega} {\bf u}_{|F} \cdot {\bf v}_{|F} dS \approx \displaystyle area(F \subset \partial \Omega) \sum_{i=1}^3 \frac{{\bf u}(R_{F,i}) 
\cdot {\bf v}(R_{F,i})}{3}, 
\end{array}
\]
\noindent where $S_{K,i}$ are the vertexes of $K$, $i=1,2,3,4$, $G_K$ is the centroid of $K$ and $R_{F,i}$ are the vertexes of a face $F \subset \partial \Omega$ of  certain tetrahedrons $K \in {\mathcal T}_h$, $i=1,2,3$.\\

Before pursuing we define the auxiliary function $\varepsilon_h$ whose value in each $K \in {\mathcal T}_h$ is constant equal to $\varepsilon(G_K)$. Furthermore 
we introduce the norms $\parallel \{\cdot\} \parallel_{\varepsilon_h,h}$ and $\parallel \{\cdot\} \parallel_{h}$ of $V_h$, given by 
$( \{\cdot\}, \{\cdot\} )_{\varepsilon_h,h}^{1/2}$ and $(\{\cdot\},\{\cdot\})_h^{1/2}$, respectively. Similarly we denote by $\parallel \{\cdot\} \parallel_{\partial \Omega,h}$ the norm defined by $(\{\cdot\},\{\cdot\})_{\partial\Omega,h}^{1,2}$. 
Then still denoting the approximation of ${\bf e}_h(k\tau)$
by ${\bf e}_h^k$, for $k=1,2,...,N$ we determine ${\bf e}_h^{k+1}$ by,
\begin{equation}\label{eq7}
  \begin{array}{l}
    \displaystyle \left( \frac{{\bf e}_h^{k+1} - 2 {\bf e}_h^k + {\bf e}_h^{k-1}}{\tau^2}, {\bf v} \right)_{\varepsilon_h,h} + (\nabla {\bf e}_h^k, \nabla {\bf v}) 
     + (\nabla \cdot \varepsilon {\bf e}_h^k, \nabla \cdot {\bf v}) - (\nabla \cdot {\bf e}_h^k, \nabla \cdot {\bf v}) \\
		 + \displaystyle \left(\frac{{\bf e}_h^{k+1} - {\bf e}_h^{k-1}}{2\tau}, {\bf v} \right)_{\partial \Omega,h} = 0 \; \forall {\bf v} \in {\bf V}_h,\\
		\\
     {\bf e}_h^0 = {\bf e}_{0h} \mbox{ and } {\bf e}_h^1 = {\bf e}_h^0 + \tau {{\bf e}_1}_h \mbox{ in } \Omega.
  \end{array}
\end{equation}

Now we recall a result given in Lemma 3 of \cite{JCAM2010}, which allows us to assert that the norm $\parallel v \parallel_{\varepsilon_h}$ is bounded above by 
$\parallel v \parallel_{\varepsilon_h,h}$ $\forall v \in V_h$. In order to prove such a result we use the barycentric coordinates $\lambda_{K,i}$ of 
tetrahedron $K \in {\mathcal T}_h$, $i=1,2,3,4$. We have $u_{|K} = \displaystyle \sum_{i=1}^4 u(S_{K,i}) \lambda_{K,i}$. Since 
$\int_k \lambda_{K,i} \lambda_{K,j} d{\bf x} = volume(K) (1+ \delta_{i,j})/20$, after straightforward manipulations we obtain,    
\[ \displaystyle \int_K \varepsilon_h u_{|K}^2 d{\bf x} \leq \displaystyle \frac{\varepsilon(G_K) volume(K)}{4}  \sum_{i=1}^4 u^2(S_{K,i}). \]
This immediately implies that 
\begin{equation}
\label{upperbound}
\parallel v \parallel_{\varepsilon_h} \leq \parallel v \parallel_{\varepsilon_h,h} \; \forall v \in V_h.
\end{equation}
For the same reason we have, %that if $\lambda_{T,i}$ corresponds to vertex $R_{T,i} \in \partial \Omega$ for a tetrahedron $T$ having a face $F$ contained in $\partial \Omega$, we have $\int_F \lambda_{T,i}^2 dS= area(F)/6$ for $i=1,2,3$.  
%Thus, similarly to \eqref{upperbound}, we have
\begin{equation}
\label{upperboundary}
\parallel v \parallel_{\partial \Omega} \leq \parallel v \parallel_{\partial \Omega,h} \; \forall v \in V_h.
\end{equation} 
%It is noteworthy that $|(\varepsilon_h u,v)_h| \leq \parallel \varepsilon_h u \parallel_h \parallel v \parallel_h, \; \forall u,v \in V_h$.

\section{Stability analysis}

In order to conveniently prepare the subsequent steps of the reliability study of scheme \eqref{eq7}, following a technique thoroughly 
exploited in \cite{Ruas}, we carry out the stability analysis of a more general form thereof, namely:
\begin{equation}\label{eq7bis}
  \begin{array}{l}
    \displaystyle \left(\frac{{\bf e}_h^{k+1} - 2 {\bf e}_h^k + {\bf e}_h^{k-1}}{\tau^2}, {\bf v} \right)_{\varepsilon_h,h} +
    (\nabla {\bf e}_h^k, \nabla {\bf v}) + (\nabla \cdot \varepsilon {\bf e}_h^k, \nabla \cdot {\bf v}) 
    - (\nabla \cdot {\bf e}_h^k, \nabla \cdot {\bf v}) \\
   + \displaystyle \left(\frac{{\bf e}_h^{k+1} - {\bf e}_h^{k-1}}{2\tau}, {\bf v} \right)_{\partial \Omega,h}= {\bf F}^k({\bf v} ) + (d^k,\nabla \cdot {\bf v}) + 
	{\bf G}^k({\bf v}) ~ \forall {\bf v} \in {\bf V}_h,\\
    \\
    {\bf e}_h^0 = {\bf e}_{0h} \mbox{ and } {\bf e}_h^1 = {\bf e}_h^0 + \tau {\bf e}_{1h} \mbox{ in } \Omega.
  \end{array}
\end{equation}
\noindent where for every $k \in \{1,2,\ldots,N-1\}$, ${\bf F}^k$ and ${\bf G}^k$ are given bounded linear functionals over ${\bf V}_h$ and the space of traces over $\partial \Omega$ of fields in ${\bf V}_h$ equipped with the norms $\| \cdot \|_h$ and $\| \cdot \|_{\partial \Omega,h}$ respectively. We denote by $| {\bf F}^k |_h$ and $| {\bf G}^k |_{\partial \Omega,h}$ the underlying norms of both functionals. $d^k$ in turn is a given function in $L^2(\Omega)$ for $k \in \{1,2,\ldots,N-1\}.$\\

Taking ${\bf v}= {\bf e}_h^{k+1} - {\bf e}_h^{k-1}$ in \eqref{eq7bis} we get for $k=1,2,\ldots,N-1$,
\begin{equation}\label{eq8}
  \begin{array}{l}
     \displaystyle \left( \frac{{\bf e}_h^{k+1} -  2{\bf e}_h^k + {\bf e}_h^{k-1}}{\tau}, \displaystyle 
    \frac{{\bf e}_h^{k+1} - {\bf e}_h^{k-1}}{\tau} \right)_{\varepsilon_h,h} + 
    (\nabla {\bf e}_h^k, \nabla {\bf e}_h^{k+1} - \nabla {\bf e}_h^{k-1}) \\
		\\
    + (\nabla \cdot \{\varepsilon - 1\} {\bf e}_h^k, \nabla \cdot {\bf e}_h^{k+1} - \nabla \cdot {\bf e}_h^{k-1})  
   + \displaystyle \left(\frac{{\bf e}_h^{k+1} - {\bf e}_h^{k-1}}{2\tau},  {\bf e}_h^{k+1} - {\bf e}_h^{k-1} \right)_{\partial \Omega.h} \\
	\\
	= {\bf F}^k({\bf e}_h^{k+1} - {\bf e}_h^{k-1}) + (d^k,\nabla \cdot \{{\bf e}_h^{k+1} - {\bf e}_h^{k-1}\}) + 
	{\bf G}^k({\bf e}_h^{k+1} - {\bf e}_h^{k-1}) 
  \end{array}
\end{equation}

Noting that ${\bf e}_h^{k+1} -  2{\bf e}_h^k + {\bf e}_h^{k-1}= ({\bf e}_h^{k+1} - {\bf e}_h^k) -({\bf e}_h^k-{\bf e}_h^{k-1})$ 
and that $  {\bf e}_h^{k+1}-{\bf e}_h^{k-1} = ({\bf e}_h^{k+1} - {\bf e}_h^k) + ({\bf e}_h^k-{\bf e}_h^{k-1})$, 
the following estimate trivially holds for equation \eqref{eq7bis}:
\begin{equation}\label{eq8bis}
  \begin{array}{l}
    \displaystyle \left\|\frac{{\bf e}_h^{k+1} - {\bf e}_h^{k}}{\tau} \right \|_{\varepsilon_h,h}^2 -
    \displaystyle \left\|\frac{{\bf e}_h^{k} - {\bf e}_h^{k-1}}{\tau} \right \|_{\varepsilon_h,h}^2
    +I_1 + I_2 + 2 \tau \displaystyle \left\| \frac{{\bf e}_h^{k+1} - {\bf e}_h^{k-1}}{2\tau} \right\|_{\partial \Omega,h}^2 \\
		\leq | {\bf F}^k |_h \left\| {\bf e}_h^{k+1} - {\bf e}_h^{k-1} \right\|_{h}  + \|  d^k\| \|\nabla \cdot ( {\bf e}_h^{k+1} - {\bf e}_h^{k-1} )\| \\
		+ | {\bf G}^k |_{\partial \Omega,h} \;\left\| {\bf e}_h^{k+1} - {\bf e}_h^{k-1} \right\|_{\partial \Omega,h} \\
		\mbox{where} \\
		I_1:= (\nabla {\bf e}_h^k, \nabla \{{\bf e}_h^{k+1} - {\bf e}_h^{k-1} \}) ; \\
		\mbox{and} \\
		I_2:= (\nabla \cdot \{\varepsilon - 1\} {\bf e}_h^k, \nabla \cdot \{{\bf e}_h^{k+1} - {\bf e}_h^{k-1}\}).
  \end{array}
\end{equation}
Next we estimate the terms $I_1$ and $I_2$ given by \eqref{eq8bis}. \\
First of all it is easy to see that
\begin{equation}\label{eq9}
  \begin{array}{l}
  I_1 = \displaystyle \frac{1}{2}  (\| \nabla {\bf e}_h^{k+1}\|^2 + \| \nabla {\bf e}_h^{k}\|^2 - \| \nabla ({\bf e}_h^{k+1} - {\bf e}_h^k)\|^2 ) \\ 
- \displaystyle \frac{1}{2} (\| \nabla {\bf e}_h^{k}\|^2 + \| \nabla {\bf e}_h^{k-1}\|^2 - \| \nabla ({\bf e}_h^{k} - {\bf e}_h^{k-1})\|^2 ). 
\end{array}
\end{equation}
Next we note that,
\begin{equation}
\label{eq100}
  \begin{array}{l}
    I_2= J_{1} + J_{2} \\
			\mbox{where } \\
 J_1:= (\nabla \cdot {\bf e}_h^k, \nabla \cdot \{{\bf e}_h^{k+1} - {\bf e}_h^{k-1}\})_{\varepsilon -1} \\
\mbox{and } \\
 J_2 := (\nabla \varepsilon \cdot {\bf e}_h^k, \nabla \cdot \{{\bf e}_h^{k+1} - {\bf e}_h^{k-1}\}).
\end{array}
\end{equation}
Similarly to \eqref{eq9} we can write,
\begin{equation}
\label{eq10}
\begin{array}{l}
 J_1 = \displaystyle \frac{1}{2}  (\| \nabla \cdot {\bf e}_h^{k+1}\|_{\varepsilon-1}^2 +  \| \nabla \cdot {\bf e}_h^{k}\|_{\varepsilon-1}^2 -
    \| \nabla \cdot({\bf e}_h^{k+1} - {\bf e}_h^k)\|_{\varepsilon-1}^2 ) \\  
  - \displaystyle \frac{1}{2}  (\| \nabla \cdot {\bf e}_h^{k}\|_{\varepsilon-1}^2 + \| \nabla \cdot {\bf e}_h^{k-1}\|_{\varepsilon-1}^2 -
    \| \nabla \cdot ({\bf e}_h^{k} - {\bf e}_h^{k-1})\|_{\varepsilon-1}^2 )
\end{array}
\end{equation}
%where $\| (\cdot) \|_{\varepsilon} = \|\sqrt{\varepsilon-1} (\cdot) \|$.\\
Now observing that $\nabla \varepsilon \equiv {\bf 0}$ on $\partial \Omega$, we integrate by parts $J_{2}$ given by \eqref{eq100}, to get
\begin{equation}
\label{eq11}
   \begin{array}{l}
     J_2 = -(\nabla \{\nabla \varepsilon \cdot {\bf e}_h^k\}, {\bf e}_h^{k+1} - {\bf e}_h^{k-1}).
\end{array}
\end{equation}
Let us rewrite $J_2$ as, 
\begin{equation}
\label{eq12bis}
   \begin{array}{l}
     J_2 = M_1 + M_2 \\
\mbox{where } \\
M_1 :=- \left (\nabla \nabla \varepsilon {\bf e}_h^k,{\bf e}_h^{k+1}- {\bf e}_h^{k-1} \right) \\
\mbox{and} \\
M_2 := - \left (\{\nabla {\bf e}_h^k\}^{\cal T} \nabla \varepsilon,{\bf e}_h^{k+1}-{\bf e}_h^{k-1} \right)
\end{array}
\end{equation}
$M_1$ in turn can be rewritten as follows:
\begin{equation}\label{eq12}
   \begin{array}{l}
    % -(\nabla \nabla \varepsilon {\bf e}_h^k, {\bf e}_h^{k+1} - {\bf e}_h^{k-1})
     M_1 = N_1 + N_2 \\
		\mbox{where} \\
		N_1:= -\tau \displaystyle \left(\nabla \nabla \varepsilon {\bf e}_h^k,\frac{ {\bf e}_h^{k+1}
     - {\bf e}_h^{k}}{\tau} \right) \\
		\mbox{and}\\
		N_2:= -\tau \displaystyle \left (\nabla \nabla \varepsilon {\bf e}_h^k,\frac{ {\bf e}_h^{k} - {\bf e}_h^{k-1}}{\tau} \right).
\end{array}
\end{equation}
Then we further observe that
\begin{equation}\label{eq13}
   %\begin{array}{l}
     N_1 = %-\tau \displaystyle \left(\nabla \nabla \varepsilon {\bf e}_h^k,\frac{ {\bf e}_h^{k+1} - {\bf e}_h^{k}}{\tau} \right) = \\
		-\tau^2 \displaystyle \sum_{i=1}^k \left (\nabla \nabla \varepsilon \frac{{\bf e}_h^i - {\bf e}_h^{i-1}}{\tau},  \frac{{\bf e}_h^{k+1} - {\bf e}_h^{k}}{\tau}       \right) - \tau \displaystyle \left(\nabla \nabla \varepsilon {\bf e}_h^0,\frac{ {\bf e}_h^{k+1} - {\bf e}_h^{k}}{\tau} \right).
%\end{array}
\end{equation}
and hence,
%and applying Young's inequality for estimation of terms $ab \leq\frac{1}{2} a ^2 + \frac{1}{2} b^2 $ 
\begin{equation*}
%      \left|(\nabla \nabla \varepsilon {\bf e}_h^k,{\bf e}_h^{k+1} - {\bf e}_h^{k}) \right| = \displaystyle \left|\tau \left(\nabla \nabla \varepsilon {\bf e}_h^k,\frac{ {\bf e}_h^{k+1} - {\bf e}_h^{k}}{\tau} \right) \right|  \\
		N_1 \geq - \tau^2 | \varepsilon |_{2,\infty} \displaystyle \left \|\frac{{\bf e}_h^{k+1} - {\bf e}_h^{k}}{\tau} \right \| \displaystyle \sum_{i=1}^k  \left \| \frac{{\bf e}_h^{i} - {\bf e}_h^{i-1}}{\tau} \right \| + \tau | \varepsilon |_{2,\infty} \| {\bf e}_h^0\|
    \left  \|\frac{{\bf e}_h^{k+1} - {\bf e}_h^{k}}{\tau} \right \|  
\end{equation*}
or yet,
\begin{equation*}
  N_1  \geq -\tau | \varepsilon |_{2,\infty} \displaystyle \left \|\frac{{\bf e}_h^{k+1} - {\bf e}_h^{k}}{\tau} \right \| \left\{ \tau \sqrt{k}  \left(  \displaystyle \sum_{i=1}^k  \left \| \frac{{\bf e}_h^{i} - {\bf e}_h^{i-1}}{\tau} \right \|^2 \right)^{1/2} + \| {\bf e}_h^0\| \right\}, 
	\end{equation*}
and noting that $k \leq T/\tau$ we get  
\begin{equation}
\label{eq14}
 \begin{array}{l}
 N_1  \geq -\tau | \varepsilon |_{2,\infty} \displaystyle \left \|\frac{{\bf e}_h^{k+1} - {\bf e}_h^{k}}{\tau} \right \| \left\{ \tau \sqrt{\frac{T}{\tau}}  \left(  \displaystyle \sum_{i=1}^k  \left \| \frac{{\bf e}_h^{i} - {\bf e}_h^{i-1}}{\tau} \right \|^2 \right)^{1/2} + \| {\bf e}_h^0\| \right\}. 
\end{array}
\end{equation}
Applying to \eqref{eq14} Young's inequality $ab \leq \delta a^2/2 + b^2/(2 \delta)$ $\forall a, b \in \Re$ and $\delta >0$ with $\delta = 1$, we easily conclude that
\begin{equation}\label{eq15}
   \begin{array}{l}
     N_1 \geq 
     - \displaystyle \frac{\tau}{2} | \varepsilon |_{2,\infty} \displaystyle  \left( 2 \left \|\frac{{\bf e}_h^{k+1} - {\bf e}_h^{k}}{\tau} \right \|^2  +  
		\tau T\displaystyle \sum_{i=1}^k  \left \| \frac{{\bf e}_h^{i} - {\bf e}_h^{i-1}}{\tau} \right\|^2 + \| {\bf e}_h^0\|^2 \right).
\end{array}
\end{equation}
Similarly to \eqref{eq15},
\begin{equation}\label{eq15bis}
   \begin{array}{l}
     N_2 \geq  
		- \displaystyle \frac{\tau}{2} | \varepsilon |_{2,\infty} \displaystyle \left( 2 \left\|\frac{{\bf e}_h^{k} - {\bf e}_h^{k-1}}{\tau} \right \|^2 + \tau T \displaystyle \sum_{i=1}^k  \left \| \frac{{\bf e}_h^{i} - {\bf e}_h^{i-1}}{\tau} \right\|^2 + \| {\bf e}_h^0\|^2 \right).
\end{array}
\end{equation}
Combining \eqref{eq15} and \eqref{eq15bis} we come up with
\begin{equation}
\label{eq15ter}
   \begin{array}{l}
     %-(\nabla \nabla \varepsilon {\bf e}_h^k,{\bf e}_h^{k+1} - {\bf e}_h^{k-1}) 
		M_1 \geq -\tau | \varepsilon |_{2,\infty} 
        \left( \displaystyle \left \|\frac{{\bf e}_h^{k+1} - {\bf e}_h^{k}}{\tau} \right \|^2 
		+\displaystyle  \left\|\frac{{\bf e}_h^{k} - {\bf e}_h^{k-1}}{\tau} \right \|^2 +  L_k \right),\\
		\mbox{where} \\
		L_k := \tau T \displaystyle \sum_{i=1}^k  \left \| \frac{{\bf e}_h^{i} - {\bf e}_h^{i-1}}{\tau} \right\|^2 + \| {\bf e}_h^0\|^2.
\end{array}
\end{equation}
As for $M_2$ given by (\ref{eq12bis}) we have:
\begin{equation*}
 % \begin{array}{l}
	%- \left (\{\nabla {\bf e}_h^k\}^{\cal T} \nabla \varepsilon ,{\bf e}_h^{k+1}-{\bf e}_h^{k-1} \right) 
	M_2 \geq - \left| \varepsilon \right|_{1,\infty} \left\| \nabla {\bf e}_h^k \right\| \left\| {\bf e}_h^{k+1} - {\bf e}_h^{k-1} \right\| \geq 
		-\tau | \varepsilon|_{1,\infty} \left\| \nabla {\bf e}_h^k \right\|
   \displaystyle \left( \left \|\frac{{\bf e}_h^{k+1} - {\bf e}_h^{k}}{\tau} \right \| + \displaystyle \left \|\frac{{\bf e}_h^{k} - {\bf e}_h^{k-1}}{\tau} \right \| \right)
\end{equation*}
or yet
\begin{equation}
\label{eq16}
 M_2 \geq 
    - \displaystyle \frac{\tau}{2}  | \varepsilon|_{1,\infty} \left( 2 \| \nabla {\bf e}_h^k \|^2 +  \displaystyle \left \|\frac{{\bf e}_h^{k+1} - {\bf e}_h^{k}}{\tau} \right \|^2  + \displaystyle \left \|\frac{{\bf e}_h^{k} - {\bf e}_h^{k-1}}{\tau} \right \|^2   \right).
\end{equation}
Now we recall \eqref{eq8bis} together with \eqref{upperboundary} and notice that for every square-integrable field ${\bf A}$ in $\Omega$ we have 
$\| {\bf A} \| \leq \| {\bf A} \|_{\varepsilon_h}$, Then taking into account that $\| \nabla \cdot {\bf v} \| \leq \sqrt{3} \| \nabla {\bf v} \|$ $\forall {\bf v} \in [H^1(\Omega)]^3$, and using Young's inequality with $\delta=\tau$, $\delta = 1/\tau$ and $\delta = \tau$, respectively, we easily 
obtain the following estimates:
%\textbf{(The second estimate is not clear: why we divide to $\tau$ here?
%  It is supposed that we should use Young's inequality for estimation of terms $ab \leq \frac{1}{2} a ^2 + \frac{1}{2} b^2 $ or?  In  the last one estimate $1/2$ also disappears. More details should be written here.  The last boundary term in  the last inequality disappears in the final estimate. Not clear, should be explained.)}
\begin{equation}
\label{eq17}
  \begin{array}{l}
    | {\bf F}^k |_h \; \| {\bf e}_h^{k+1} - {\bf e}_h^{k-1}\|  \leq \displaystyle \frac{\tau}{2} | {\bf F}^k |^2_h 
     + \displaystyle \tau \left( \left \|\frac{{\bf e}_h^{k+1} - {\bf e}_h^{k}}{\tau} \right \|^2_h  + \displaystyle \left \|\frac{{\bf e}_h^{k} - {\bf e}_h^{k-1}}{\tau} \right \|^2_h \right),\\
		\\
    \| d^k\| \; \|\nabla \cdot ( {\bf e}_h^{k+1} - {\bf e}_h^{k-1} ) \| %\leq  \sqrt{3}  \| d^k\| \; \|\nabla ( {\bf e}_h^{k+1} - {\bf e}_h^{k-1} ) \|
    %+ \frac{\tau}{2} \| \nabla ({\bf e}_h^{k+1} - {\bf e}_h^{k-1})\|^2 
		\leq \displaystyle \frac{3}{2 \tau} \| d^k\|^2 + 
		\tau  ( \| \nabla {\bf e}_h^{k+1} \|^2 + \| \nabla {\bf e}_h^{k-1})\|^2), \\
    \\
		|{\bf G}^k |_{\partial \Omega,h} \displaystyle \left\| {\bf e}_h^{k+1} - {\bf e}_h^{k-1} \right\|_{\partial \Omega,h} \leq  
    \displaystyle \frac{\tau}{2} |{\bf G}^k |_{\partial \Omega,h}^2 
		+ 2 \tau \displaystyle \left\| \frac{{\bf e}_h^{k+1} - {\bf e}_h^{k-1}}{2 \tau} \right\|_{\partial \Omega,h}^2.\,
  \end{array}
\end{equation}
where in the first and the second inequality we also used the fact that $\| {\bf A} \pm {\bf B} \|^2 \leq 2 (\| {\bf A} \|^2 + \| {\bf B} \|^2 )$ for all square-integrable fields 
${\bf A}$ and ${\bf B}$.\\
%Let us definitively take $\tau$ such that the condition \eqref{CFL} holds.\\
Now we collect \eqref{eq9}, \eqref{eq100}, \eqref{eq10}, \eqref{eq11}, \eqref{eq12bis} and \eqref{eq15ter}, \eqref{eq16}, \eqref{eq17} to plug them into \eqref{eq8bis}. Using the fact that $\| \cdot \| \leq \| \cdot \|_h \leq\| \cdot \|_{\varepsilon_h}\leq \| \cdot \|_{\varepsilon_h,h}$  
%summing up from $k=1$ through $k=m$ and taking into account \eqref{eq18}, 
we obtain for $1 \leq k \leq N-1$:
%Equation \eqref{eq19} can be rewritten as follows, for $2 \leq m \leq N-1$:
\begin{equation}
\label{eq19} 
  \begin{array}{l}
   ( A_k - B_k ) + ( C_k - D_k ) 
	\leq \displaystyle \frac{\tau}{2}| {\bf F}^k|^2_h +
 \displaystyle \frac{3}{2 \tau} \| d^k\|^2  + \displaystyle \frac{\tau}{2} \left| {\bf G}^k \right|^2_{\partial \Omega_h}, \\
%+ T | \varepsilon |_{2,\infty} \| {\bf e}_h^0\|^2 \\
  \mbox{where} \\
	A_k:=  \displaystyle \left\| \frac{{\bf e}_h^{k+1} - {\bf e}_h^{k}}{\tau} \right\|^2_{\varepsilon_h,h} 
	- \displaystyle \frac{\tau}{2} \left( 2 + | \varepsilon |_{1,\infty} +
  2 | \varepsilon |_{2,\infty} \right) \displaystyle \left\| \frac{{\bf e}_h^{k+1} - {\bf e}_h^{k}}{\tau} \right\|^2_{\varepsilon_h,h} \\
	+ \displaystyle \frac{1}{2} \left\{ (1-2\tau)\|\nabla  {\bf e}_h^{k+1}\|^2  + (1-\tau | \varepsilon |_{1,\infty}) \| \nabla {\bf e}_h^k \|^2 \right\} \\
	\\ 
	B_k := \displaystyle  \left\| \frac{{\bf e}_h^{k} - {\bf e}_h^{k-1}}{\tau}\right \|^2_{\varepsilon_h,h} + \displaystyle \frac{\tau}{2} 
	\left( 2 + | \varepsilon |_{1,\infty} + 2| \varepsilon|_{2,\infty} \right) \displaystyle \left\| \frac{{\bf e}_h^{k} - {\bf e}_h^{k-1}}{\tau}\right \|^2_{\varepsilon_h,h} \\
	+ \displaystyle \frac{1}{2} \left\{ \left(1 + 2\tau \right) \| \nabla {\bf e}_h^{k-1}\|^2 + \left( 1 + \tau |\varepsilon|_{1,\infty} \right) \| \nabla {\bf e}_h^k\|^2 \right\} + \tau | \varepsilon |_{2,\infty} L_k\\
	\\
	C_k := \displaystyle -\frac{1}{2} \left\{ 
  \|\nabla ({\bf e}_h^{k+1} - {\bf e}_h^k)\|^2 + \| \nabla \cdot ({\bf e}_h^{k+1} - {\bf e}_h^k) \|_{\varepsilon-1}^2 \right\} +  \displaystyle 
	\frac{1}{2} \left( \| \nabla \cdot {\bf e}_h^{k+1}\|_{\varepsilon-1}^2 + \| \nabla \cdot {\bf e}_h^k\|_{\varepsilon-1}^2 \right) \\
	\\
	D_k := \displaystyle -\frac{1}{2} \left\{ 
  \|\nabla ({\bf e}_h^k - {\bf e}_h^{k-1})\|^2 + \| \nabla \cdot ({\bf e}_h^k - {\bf e}_h^{k-1})\|_{\varepsilon-1}^2 \right\} + \displaystyle \frac{1}{2} 
  \left( \| \nabla \cdot {\bf e}_h^k\|_{\varepsilon-1}^2 + \| \nabla \cdot {\bf e}_h^{k-1} \|_{\varepsilon-1}^2 \right).
\end{array}
\end{equation}
Setting 
\begin{equation}
\label{eta}
\begin{array}{l}
\eta :=  2 + | \varepsilon |_{1,\infty} + 2| \varepsilon|_{2,\infty}; \\
\theta := | \varepsilon |_{1,\infty}; \\
\rho := T^2 | \varepsilon |_{2,\infty},
\end{array}
\end{equation}
we can rewrite $A_k$ and $B_k$ as follows: 
\begin{equation}
\label{eq20}
  \begin{array}{l}
	A_k:=  \displaystyle \left\| \frac{{\bf e}_h^{k+1} - {\bf e}_h^{k}}{\tau} \right\|^2_{\varepsilon_h,h} 
	- \displaystyle \frac{\tau}{2} \eta \displaystyle \left\| \frac{{\bf e}_h^{k+1} - {\bf e}_h^{k}}{\tau} \right\|^2_{\varepsilon_h,h} \\
	+ \displaystyle \frac{1}{2} \left(\|\nabla  {\bf e}_h^{k+1}\|^2  + \| \nabla {\bf e}_h^k \|^2 \right)  
	-\displaystyle  \tau  \left( \|\nabla  {\bf e}_h^{k+1}\|^2 + \frac{\theta}{2} \| \nabla {\bf e}_h^k \|^2 \right); \\
	\\
	B_k:=  \displaystyle \left\| \frac{{\bf e}_h^{k} - {\bf e}_h^{k-1}}{\tau} \right\|^2_{\varepsilon_h,h} 
	+ \displaystyle \frac{\tau}{2} \eta \displaystyle \left\| \frac{{\bf e}_h^{k} - {\bf e}_h^{k-1}}{\tau} \right\|^2_{\varepsilon_h,h} \\
	 %\left\| \frac{{\bf e}_h^{k} - {\bf e}_h^{k-1}}{\tau} \right\|^2_{\varepsilon_h,h} \\
	+ \displaystyle \frac{1}{2} \left(\|\nabla  {\bf e}_h^{k-1}\|^2  + \| \nabla {\bf e}_h^k \|^2 \right)  
	+\displaystyle  \tau  \left( \|\nabla  {\bf e}_h^{k-1}\|^2 + \frac{\theta}{2} \| \nabla {\bf e}_h^k \|^2 \right) + \displaystyle \frac{\tau \rho}{T^2} L_k \\
\end{array}
\end{equation}
Then we note that for $1 \leq k \leq m - 1$ with $2 \leq m \leq N-1$,
\begin{equation}
\label{eq21}
   \begin{array}{l}
	A_{k-1} - B_{k} = - \displaystyle \tau \eta \left\| \frac{{\bf e}_h^{k} - {\bf e}_h^{k-1}}{\tau} \right\|^2_{\varepsilon_h,h} 
	-\displaystyle  \tau \left( 1+\frac{\theta}{2} \right) \left( \|\nabla  {\bf e}_h^{k}\|^2 + \| \nabla {\bf e}_h^{k-1} \|^2 \right) - \displaystyle \frac{\tau \rho}{T^2} L_k.
	\end{array}
\end{equation}
	It follows that
	\begin{equation}
	\label{eq21bis}
	\begin{array}{l}
	 \displaystyle \sum_{k=1}^m (A_k-B_k)  = A_m - B_1 + \displaystyle \sum_{k=2}^{m} (A_{k-1}-B_{k})  \\
	= \displaystyle \frac{1}{2} \left( 2- \tau \eta \right) \displaystyle \left\| \frac{{\bf e}_h^{m+1} - {\bf e}_h^{m}}{\tau} \right\|^2_{\varepsilon_h,h} 
   + \displaystyle  \frac{1}{2} \left( 1 - 2 \tau \right) \|\nabla  {\bf e}_h^{m+1}\|^2 + 
	\displaystyle \frac{1}{2} \left( 1 - \tau \theta \right) \| \nabla {\bf e}_h^m \|^2 \\ 
	- \displaystyle \frac{1}{2} \left( 2 + \tau \eta \right) \left\| \frac{{\bf e}_h^{1} - {\bf e}_h^{0}}{\tau} \right\|^2_{\varepsilon_h,h} 
	- \displaystyle \frac{1}{2} \left( 1 + 2 \tau \right) \|\nabla  {\bf e}_h^{0}\|^2  - \frac{1}{2} \left( 1 +\tau \theta \right)  
	\| \nabla {\bf e}_h^1 \|^2 - \displaystyle \frac{\tau \rho}{T^2} L_1 \\  
	- \displaystyle \sum_{k=2}^m  \left\{ \displaystyle \tau \eta \left\| \frac{{\bf e}_h^{k} - {\bf e}_h^{k-1}}{\tau} \right\|^2_{\varepsilon_h,h} + \displaystyle \frac{\tau}{2} \left( 2+ \theta \right) \left( \|\nabla  {\bf e}_h^{k}\|^2 + \| \nabla {\bf e}_h^{k-1} \|^2 \right) + \displaystyle \frac{\tau \rho}{T^2} L_k \right\}.
	\end{array}
	\end{equation}
	Now we extend to $k=1$ the summation range on the right hand side of \eqref{eq21bis} to obtain,
	\begin{equation}
	\label{eq21ter}
	  \begin{array}{l}
	 \displaystyle \sum_{k=1}^m (A_k-B_k) \\
	= \displaystyle \frac{1}{2} \left( 2- \tau \eta \right) \displaystyle \left\| \frac{{\bf e}_h^{m+1} - {\bf e}_h^{m}}{\tau} \right\|^2_{\varepsilon_h,h} 
    + \displaystyle \frac{1}{2} \left( 1 - 2 \tau \right) \|\nabla  {\bf e}_h^{m+1}\|^2 + 
	\displaystyle \frac{1}{2} \left( 1 - \tau \theta \right) \| \nabla {\bf e}_h^m \|^2 \\ 
	- \displaystyle \frac{1}{2} \left( 2 - \tau \eta \right) \left\| \frac{{\bf e}_h^{1} - {\bf e}_h^{0}}{\tau} \right\|^2_{\varepsilon_h,h} 
	- \displaystyle \frac{1}{2} \left( 1 - \tau \theta \right) \|\nabla  {\bf e}_h^{0}\|^2 - \frac{1}{2} \left( 1 - 2 \tau \right) 
	\| \nabla {\bf e}_h^1 \|^2  \\  
	- \displaystyle \sum_{k=1}^m  \left\{ \displaystyle \tau \eta \left\| \frac{{\bf e}_h^{k} - {\bf e}_h^{k-1}}{\tau} \right\|^2_{\varepsilon_h,h} + \displaystyle  \frac{\tau}{2} \left( 2+ \theta \right) \left( \|\nabla  {\bf e}_h^{k}\|^2 + \| \nabla {\bf e}_h^{k-1} \|^2 \right) + \displaystyle \frac{\tau \rho}{T^2}  L_k\right\}.
	\end{array}
	\end{equation}
	Moreover since $C_{k-1} = D_k$ for all $m \geq k \geq 2$, recalling \eqref{eq19} we easily derive,
	\begin{equation}
	\label{eq22}
  \begin{array}{l} 
	\displaystyle \sum_{k=1}^m (C_k - D_k) = C_m - D_1 \\
  = - \displaystyle \frac{1}{2} 
  \left\{ \|\nabla ({\bf e}_h^{m+1} - {\bf e}_h^m)\|^2 + \| \nabla \cdot ({\bf e}_h^{m+1} - {\bf e}_h^m)\|_{\varepsilon-1}^2 \right\} 
	+ \displaystyle \frac{1}{2} \left( \| \nabla \cdot {\bf e}_h^{m+1}\|_{\varepsilon-1}^2 + \| \nabla \cdot {\bf e}_h^m\|_{\varepsilon-1}^2 \right) \\
	+ \displaystyle \frac{1}{2} \left\{ 
  \|\nabla ({\bf e}_h^1 - {\bf e}_h^{0})\|^2 + \| \nabla \cdot ({\bf e}_h^1 - {\bf e}_h^{0})\|_{\varepsilon-1}^2 \right\}  	
	- \displaystyle \frac{1}{2} 
  \left( \| \nabla \cdot {\bf e}_h^1\|_{\varepsilon-1}^2  + \| \nabla \cdot {\bf e}_h^{0} \|_{\varepsilon-1}^2 \right). 
\end{array}
\end{equation}
Combining \eqref{eq21ter} and \eqref{eq22}, we obtain for $2 \leq m \leq N-1$: 
\begin{equation}
 \label{eq22bis}
   \begin{array}{l} 
	\displaystyle \sum_{k=1}^m (A_k + C_k) - \displaystyle \sum_{k=1}^m (B_k + D_k)\\ 
  = \displaystyle \frac{1}{2} \left( 2- \tau \eta \right) \displaystyle \left\| \frac{{\bf e}_h^{m+1} - {\bf e}_h^{m}}{\tau} \right\|^2_{\varepsilon_h,h} 
    + \displaystyle  \frac{1}{2} \left( 1 - 2 \tau \right) \|\nabla  {\bf e}_h^{m+1}\|^2 + 
	\displaystyle  \frac{1}{2} \left( 1 - \tau \theta \right) \| \nabla {\bf e}_h^m \|^2 \\ 
	- \displaystyle \frac{1}{2} \left( 2 - \tau \eta \right) \left\| \frac{{\bf e}_h^{1} - {\bf e}_h^{0}}{\tau} \right\|^2_{\varepsilon_h,h} 
	- \displaystyle \frac{1}{2} \left( 1 - \tau \theta \right) \|\nabla  {\bf e}_h^{0}\|^2 - \frac{1}{2} \left( 1 - 2 \tau \right) 
	\| \nabla {\bf e}_h^1 \|^2  \\  
	- \displaystyle \sum_{k=1}^m  \left\{ \displaystyle \tau \eta \left\| \frac{{\bf e}_h^{k} - {\bf e}_h^{k-1}}{\tau} \right\|^2_{\varepsilon_h,h} + \displaystyle  \frac{\tau}{2} \left( 2 + \theta \right) \left( \|\nabla  {\bf e}_h^{k}\|^2 + \| \nabla {\bf e}_h^{k-1} \|^2 \right) + \displaystyle \frac{\tau \rho}{T^2} L_k \right\} \\
	- \displaystyle \frac{1}{2} 
  \left\{ \|\nabla ({\bf e}_h^{m+1} - {\bf e}_h^m)\|^2 + \| \nabla \cdot ({\bf e}_h^{m+1} - {\bf e}_h^m)\|_{\varepsilon-1}^2 \right\} 
	+ \displaystyle \frac{1}{2} \left( \| \nabla \cdot {\bf e}_h^{m+1}\|_{\varepsilon-1}^2 + \| \nabla \cdot {\bf e}_h^m\|_{\varepsilon-1}^2 \right) \\
	+ \displaystyle \frac{1}{2} \left\{ 
  \|\nabla ({\bf e}_h^1 - {\bf e}_h^{0})\|^2 + \| \nabla \cdot ({\bf e}_h^1 - {\bf e}_h^{0})\|_{\varepsilon-1}^2 \right\}  	
	- \displaystyle \frac{1}{2} 
  \left( \| \nabla \cdot {\bf e}_h^1\|_{\varepsilon-1}^2  + \| \nabla \cdot {\bf e}_h^{0} \|_{\varepsilon-1}^2 \right).
	\end{array}
\end{equation}
On the other hand, recalling \eqref{eq15ter} we note that for $2 \leq m \leq N-1$ 
\begin{equation}
\label{eq180}
\begin{array}{l}
\displaystyle \sum_{k=1}^m  L_k 
\leq \displaystyle \sum_{k=1}^m  \left( \tau T \displaystyle \sum_{i=1}^m \left \|\frac{{\bf e}_h^i - {\bf e}_h^{i-1}}{\tau} \right \|^2 + \| {\bf e}_h^0\|^2 \right)  
= T m \tau \displaystyle \sum_{k=1}^m  \left \|\frac{{\bf e}_h^k - {\bf e}_h^{k-1}}{\tau} \right \|^2 + m \| {\bf e}_h^0\|^2 
%\leq T^2 \displaystyle \sum_{k=1}^m  \left \|\frac{{\bf e}_h^k - {\bf e}_h^{k-1}}{\tau} \right \|^2 + \frac{T}{\tau} \| {\bf e}_h^0\|^2.
\end{array}
\end{equation}
In short since $m \tau \leq T$, from \eqref{eq180} we easily derive for $2 \leq m \leq N-1$: 
\begin{equation}
\label{eq18}
\begin{array}{l}
\displaystyle \frac{\tau \rho}{T^2} \sum_{k=1}^m  L_k \leq 
 \tau \rho \displaystyle \sum_{k=1}^m  \left \|\frac{{\bf e}_h^k - {\bf e}_h^{k-1}}{\tau} \right \|^2 + \frac{\rho}{T} \| {\bf e}_h^0\|^2.
\end{array}
\end{equation}
Plugging \eqref{eq18} into \eqref{eq22bis} and summing up both sides of \eqref{eq19} from $k=1$ through $k=m$ for $2 \leq m \leq N-1$ by using \eqref{eq20}, 
\eqref{eq22bis} yields: 

\begin{equation}
\label{eq22ter}
  \begin{array}{l} 
	 \displaystyle \frac{1}{2} \left( 2- \tau \eta \right) \displaystyle \left\| \frac{{\bf e}_h^{m+1} - {\bf e}_h^{m}}{\tau} \right\|^2_{\varepsilon_h,h} 
    + \displaystyle \frac{1}{2} \left( 1 - 2 \tau \right) \|\nabla  {\bf e}_h^{m+1}\|^2 + 
	\displaystyle \frac{1}{2} \left( 1 - \tau \theta \right) \| \nabla {\bf e}_h^m \|^2 \\ 
	- \displaystyle \frac{1}{2} \left( 2 - \tau \eta \right) \left\| \frac{{\bf e}_h^{1} - {\bf e}_h^{0}}{\tau} \right\|^2_{\varepsilon_h,h} 
	- \displaystyle \frac{1}{2} \left( 1 - \tau \theta \right) \|\nabla  {\bf e}_h^{0}\|^2 - \frac{1}{2} \left( 1 - 2 \tau \right) 
	\| \nabla {\bf e}_h^1 \|^2  \\  
	- \displaystyle \sum_{k=1}^m  \left\{ \displaystyle \tau \left( \eta +  \rho \right) \left\| \frac{{\bf e}_h^{k} - {\bf e}_h^{k-1}}{\tau} \right\|^2_{\varepsilon_h,h} + \displaystyle \frac{\tau}{2} \left( 2 + \theta \right) \left( \|\nabla  {\bf e}_h^{k}\|^2 + \| \nabla {\bf e}_h^{k-1} \|^2 \right) \right\} \\
	- \displaystyle \frac{1}{2} 
  \left\{ \|\nabla ({\bf e}_h^{m+1} - {\bf e}_h^m)\|^2 + \| \nabla \cdot ({\bf e}_h^{m+1} - {\bf e}_h^m)\|_{\varepsilon-1}^2 \right\} 
	+ \displaystyle \frac{1}{2} \left( \| \nabla \cdot {\bf e}_h^{m+1}\|_{\varepsilon-1}^2 + \| \nabla \cdot {\bf e}_h^m\|_{\varepsilon-1}^2 \right) \\
	+ \displaystyle \frac{1}{2} \left\{ 
  \|\nabla ({\bf e}_h^1 - {\bf e}_h^{0})\|^2 + \| \nabla \cdot ({\bf e}_h^1 - {\bf e}_h^{0})\|_{\varepsilon-1}^2 \right\}  	
	- \displaystyle \frac{1}{2} 
  \left( \| \nabla \cdot {\bf e}_h^1\|_{\varepsilon-1}^2  + \| \nabla \cdot {\bf e}_h^{0} \|_{\varepsilon-1}^2 \right)	 \\
	\leq \displaystyle \sum_{k=1}^m \left(\frac{\tau}{2}| {\bf F}^k|^2_h +
 \displaystyle \frac{3}{2 \tau} \| d^k\|^2  + \displaystyle \frac{\tau}{2} \left| {\bf G}^k \right|^2_{\partial \Omega_h} \right) 
+ \displaystyle \frac{\rho}{T} \| {\bf e}_h^0\|^2.
	\end{array}
\end{equation}
Thus taking into account that ${\bf e}_h^1 - {\bf e}_h^0 = \tau {\bf e}_{1h}$, leaving on the left hand side only the terms with superscripts $m+1$ and $m$,  
and increasing the coefficients of $\| \nabla {\bf e}_h^j \|^2$ for $j=0,1$ and $\| {\bf e}_{1h} \|_{\varepsilon_h,h}^2$, we derive for $2 \leq m \leq N-1$:
\begin{equation}
\label{eq23}
  \begin{array}{l}
	\displaystyle \frac{1}{2} \left( 2 - \tau \eta\right) \displaystyle 
  \left\| \frac{{\bf e}_h^{m+1} - {\bf e}_h^{m}}{\tau} \right\|^2_{\varepsilon_h,h} 
	+ \displaystyle \frac{1}{2} \left( 1 - 2 \tau \right) \|\nabla  {\bf e}_h^{m+1}\|^2  +
  \displaystyle \frac{1}{2} \left(1 - \tau \theta \right) \| \nabla {\bf e}_h^m \|^2 \\ 
	- \displaystyle \frac{1}{2} \left\{ 
  \|\nabla ({\bf e}_h^{m+1} - {\bf e}_h^m)\|^2 + \| \nabla \cdot ({\bf e}_h^{m+1} - {\bf e}_h^m)\|_{\varepsilon-1}^2 \right\} + \displaystyle \frac{1}{2} 
  \left( \| \nabla \cdot {\bf e}_h^{m+1}\|_{\varepsilon-1}^2 + \| \nabla \cdot {\bf e}_h^m\|_{\varepsilon-1}^2 \right) \\ 
	\leq 
 \displaystyle \sum_{k=1}^m  \displaystyle \left\{ \tau \left(\eta + \rho \right) 
	\displaystyle \left\| \frac{{\bf e}^k_h-{\bf e}_h^{k-1}}{\tau} \right\|^2_{\varepsilon_h,h} 
   + \displaystyle \frac{\tau}{2}  \left( 2 + \theta \right) \left( \| \nabla {\bf e}_h^{k}\|^2 + \| \nabla {\bf e}_h^{k-1}\|^2 \right) \right. \\
 \left. + \displaystyle \frac{\tau}{2} | {\bf F}^k |^2_h + \displaystyle \frac{3}{2 \tau} \| d^k\|^2  + \displaystyle \frac{\tau}{2} \left| {\bf G}^k \right|^2_{\partial \Omega,h} \right\} 
+ \displaystyle \left\|{\bf e}_{1h} \right\|^2_{\varepsilon_h,h}  
	+ \displaystyle \frac{\tau^2}{2} \left( 
  \|\nabla {\bf e}_{1h} \|^2 + \| \nabla \cdot {\bf e}_{1h} \|_{\varepsilon-1}^2 \right) \\
		+	\displaystyle \frac{1}{2} \left( \|\nabla  {\bf e}_h^{1} \|^2 
	+ \| \nabla {\bf e}_h^0 \|^2 \right)	+ \displaystyle \frac{1}{2} \left( \| \nabla \cdot {\bf e}_h^1\|_{\varepsilon-1}^2  + \| \nabla \cdot {\bf e}_h^{0} \|_{\varepsilon-1}^2 \right) + \displaystyle \frac{\rho}{T} \| {\bf e}_h^0\|^2. 
\end{array}
\end{equation}
Now we recall a classical inverse inequality (cf. \cite{Ciarlet}), according to which,
\begin{equation}
\label{inversineq} 
\| \nabla v \| \leq C h^{-1} \| v \| \leq C h^{-1} \| v \|_{\varepsilon_h,h} \mbox{ for all } v \in V_h,
\end{equation} 
where $C$ is a mesh-independent constant, and we apply the trivial upper bound $\| \nabla \cdot {\bf v} \|_{\varepsilon-1} \leq \sqrt{3 \| \varepsilon - 1 \|_{\infty} } \| \nabla {\bf v} \|$ for all ${\bf v} \in {\bf V}_h$. \\
Let us assume that $\tau$ satisfies the following CFL-condition: 
\begin{equation} 
\label{CFL}
\tau \leq h /\nu \mbox{ with } \nu=C(1+ 3 \| \varepsilon - 1 \|_{\infty})^{1/2}.
\end{equation}
Then we have, $\forall {\bf v} \in {\bf V}_h$:
\begin{equation}
\label{eq16bis}
 \| \nabla {\bf v} \|^2 + \| \nabla \cdot {\bf v} \|_{\varepsilon-1}^2 \leq 
\displaystyle \nu^2 \displaystyle \frac{\tau^2}{h^2} \frac{\| {\bf v} \|^2}{\tau^2} \leq \displaystyle \frac{\| {\bf v} \|^2_{\varepsilon_h,h}}{\tau^2}   
\end{equation}
This means that  
\begin{equation}
\label{eq16ter}
\|\nabla ({\bf e}_h^{m+1} - {\bf e}_h^m)\|^2 + \| \nabla \cdot ({\bf e}_h^{m+1} - {\bf e}_h^m)\|_{\varepsilon-1}^2 \leq \displaystyle \left\| \frac{{\bf e}_h^{m+1} - {\bf e}_h^{m}}{\tau} \right\|^2_{\varepsilon_h,h}
\end{equation}
Plugging \eqref{eq16ter} into \eqref{eq23} we come up with,
\begin{equation}
\label{eq23bis}
  \begin{array}{l}
	 \displaystyle \frac{1}{2} \left( 1- \tau \eta \right) \displaystyle 
  \left\| \frac{{\bf e}_h^{m+1} - {\bf e}_h^{m}}{\tau} \right\|^2_{\varepsilon_h,h} 
   + \displaystyle \frac{1}{2} \left(1 - 2 \tau \right)\|\nabla  {\bf e}_h^{m+1}\|^2  +
  \displaystyle \frac{1}{2} \left(1 - \tau \theta \right) \| \nabla {\bf e}_h^m \|^2 \\ 
	\leq 
	\displaystyle \sum_{k=1}^m \left\{ \tau \left( \eta + \rho \right)  
	\displaystyle \left\| \frac{{\bf e}_h^{k} - {\bf e}_h^{k-1}}{\tau}\right \|^2_{\varepsilon_h,h} 
  + \displaystyle \frac{\tau}{2} \left( 2 + \theta \right) \left( \| \nabla {\bf e}_h^{k}\|^2 + \| \nabla {\bf e}_h^{k-1}\|^2 \right) \right.  \\
 \left. + \displaystyle \frac{\tau}{2}| {\bf F}^k|^2_h + \displaystyle \frac{3}{2 \tau} \| d^k\|^2  + \displaystyle \frac{\tau}{2} \left| {\bf G}^k \right|^2_{\partial \Omega,h} \right\}  + \displaystyle \left\|{\bf e}_{1h} \right\|^2_{\varepsilon_h,h} 
+ \displaystyle \frac{\tau^2}{2} \left( \|\nabla {\bf e}_{1h} \|^2 + \| \nabla \cdot {\bf e}_{1h} \|_{\varepsilon-1}^2 \right) \\
	+ \displaystyle \frac{1}{2} \left( \| \nabla {\bf e}_h^1 \|^2 + \|\nabla  {\bf e}_h^{0} \|^2 \right) 
	+ \displaystyle \frac{1}{2} \left( \| \nabla \cdot {\bf e}_h^1\|_{\varepsilon-1}^2  + \| \nabla \cdot {\bf e}_h^{0} \|_{\varepsilon-1}^2 \right) 
	+ \displaystyle \frac{\rho}{T} \| {\bf e}_h^0\|^2 . 
\end{array}
\end{equation}
Next we note that both $1 - 2 \tau$ and $1-\tau \theta$ are bounded below by $1- \tau \eta$ and moreover $\frac{1}{2}(2+\theta)$ is obviously bounded above by 
$\eta+\rho$. Therefore it is easy to see that \eqref{eq23bis} can be transformed into : 
\begin{equation}
\label{eq24}
  \begin{array}{l}
	 \displaystyle \frac{1}{2} \left( 1 - \tau \eta \right) \displaystyle \left( 
  \left\| \frac{{\bf e}_h^{m+1} - {\bf e}_h^{m}}{\tau} \right\|^2_{\varepsilon_h,h} 
   + \|\nabla  {\bf e}_h^{m+1}\|^2 + \| \nabla {\bf e}_h^m \|^2 \right)  \\ 
	\leq S_N + E_0 +
	\displaystyle \sum_{k=1}^m \tau \left( \eta + \rho \right) \displaystyle \left( \left\| \frac{{\bf e}_h^{k} - {\bf e}_h^{k-1}}{\tau}\right \|^2_{\varepsilon_h,h} 
  + \| \nabla {\bf e}_h^{k}\|^2 + \| \nabla {\bf e}_h^{k-1}\|^2 \right),
	\end{array}
\end{equation}
\noindent where
\begin{equation}
\label{eq24bis}
  \begin{array}{l} 
	S_N := \displaystyle \sum_{k=1}^{N-1} \left( \displaystyle \frac{\tau}{2} | {\bf F}^k|^2_h + \displaystyle \frac{3}{2 \tau} \| d^k\|^2  + \displaystyle \frac{\tau}{2} \left| {\bf G}^k \right|^2
	_{\partial \Omega,h} \right) \\
	\\
	E_0 := \displaystyle \left\|{\bf e}_{1h} \right\|^2_{\varepsilon_h,h} 
	+ \displaystyle \frac{\tau^2}{2} \left( 
  \|\nabla {\bf e}_{1h} \|^2 + \| \nabla \cdot {\bf e}_{1h} \|_{\varepsilon-1}^2 \right)  \\
		+ \displaystyle \frac{1}{2} \left( \| \nabla {\bf e}_h^1 \|^2 +  \|\nabla  {\bf e}_h^{0} \|^2 
	\right)
	+ \displaystyle \frac{1}{2} \left( \| \nabla \cdot {\bf e}_h^1\|_{\varepsilon-1}^2  + \| \nabla \cdot {\bf e}_h^{0} \|_{\varepsilon-1}^2 \right) 
	+ \displaystyle \frac{\rho}{T} \| {\bf e}_h^0\|^2 . 
\end{array}
\end{equation}
Let us assume that $\tau \leq \displaystyle \frac{1}{2 \eta}$. Then from \eqref{eq24bis} we have 
\begin{equation}
\label{eq25}
  \begin{array}{l} 
  \displaystyle \left\| \frac{{\bf e}_h^{m+1} - {\bf e}_h^{m}}{\tau} \right\|^2_{\varepsilon_h,h} 
   + \|\nabla  {\bf e}_h^{m+1}\|^2 + \| \nabla {\bf e}_h^m \|^2  \\ 
	\leq 4 \left\{ S_N +  E_0 +
	\displaystyle \sum_{k=1}^m  \tau \left( \eta + \rho \right) \displaystyle \left( \left\| \frac{{\bf e}_h^{k} - {\bf e}_h^{k-1}}{\tau}\right \|^2_{\varepsilon_h,h} 
  + \| \nabla {\bf e}_h^{k}\|^2 + \| \nabla {\bf e}_h^{k-1}\|^2 \right) \right\}. 
\end{array}
\end{equation}
Now setting 
\begin{equation}
\label{eq26}
  \beta= 4 N \tau (\eta + \rho) = 4 T \left\{ 2 + | \varepsilon |_{1,\infty} + (2 + T^2) | \varepsilon |_{2,\infty} \right\}, 
\end{equation}
from the discrete Gr\"onwall's Lemma and \eqref{upperbound} from \eqref{eq25} we derive for all $m \leq N-1$:
\begin{equation}
\label{stability}
\boxed{
  \displaystyle  \left\| \frac{{\bf e}_h^{m+1} - {\bf e}_h^m}{\tau} \right\|^2 + \| \nabla {\bf e}_h^{m+1}\|^2 + \| \nabla {\bf e}_h^{m}\|^2 \leq 4 (S_N + E_0) e^{\beta},
}
\end{equation}
\noindent as long as %\eqref{CFL} is satisfied and $\tau \leq \frac{1}{2 \eta}$, that is 
$\tau \leq \min\{h/\nu ,1/(2 \eta)\}$, where $\nu$, $\eta$-$\rho$ and $\beta$ are defined in \eqref{CFL}, \eqref{eta} and \eqref{eq26}, respectively, and in the expression of $E_0$, ${\bf e}^1_h$ is to be replaced by ${\bf e}^0_h + \tau {\bf e}_{1h}$. \rule{2mm}{2mm}

\section{Scheme consistency}
Before pursuing the reliability study of our scheme we need some approximation results related to the Maxwell's equations. The arguments employed in this section 
found their inspiration in Thom\'ee \cite{Thomee} and in Ruas \cite{Ruas}.

\subsection{Preliminaries}

Henceforth we assume that $\Omega$ is a convex polyhedron. In this case one may reasonably assume that for every ${\bf g} \in [L^2_0(\Omega)]^3$, 
where $L^2_0(\Omega)$ is the subspace of $L^2(\Omega)$ consisting of functions whose integral in $\Omega$ equals zero,  
the solution ${\bf v}_{\bf g} \in [H^1(\Omega) \cap L^2_0(\Omega)]^3$ of the equation 
%stationary counterpart of Maxwell's equations, namely
\begin{equation}
\label{vectpoisson}
\begin{array}{ll}
-\nabla^2 {\bf v}_{\bf g} - \nabla [(\varepsilon - 1) \nabla \cdot {\bf v}_{\bf g}] = {\bf g} & \mbox{ in } \Omega\\
 \partial_n {\bf v}_{\bf g} = {\bf 0} & \mbox{ on } \partial \Omega.
\end{array}
\end{equation}
belongs to $[H^2(\Omega)]^3$. \\
Another result that we take for granted in this section is the existence of a constant $C$ such that,
\begin{equation}
\label{regularity} 
\forall {\bf g} \in [L^2_0(\Omega)]^3 \mbox{ it holds } \| {\mathcal H}({\bf v}_{\bf g}) \| \leq C \| {\bf g} \| , 
\end{equation}
where %${\bf v}_{\bf g}$ is the solution of \eqref{vectpoisson} and 
${\mathcal H}(\cdot)$ is the Hessian of a function or field. \eqref{regularity} is a result whose grounds can be found in analogous inequalities applying to the scalar Poisson problem and to %the linear elasticity system, or yet 
the Stokes system. In fact \eqref{vectpoisson} can be viewed as a problem half way between the (vector) Poisson problem and a sort of generalized Stokes system, both with homogeneous Neumann boundary conditions. %Indeed let us define $\Omega_{\varepsilon} := \{{\bf x} | \; \varepsilon({\bf x}) > 1 \}$. 
In order to create such a system we replace in \eqref{vectpoisson} $- (\varepsilon - 1) \nabla \cdot {\bf v}_{\bf g}$ by 
a fictitious pressure $p_{\bf g}$. Then the resulting equation %by $p_{\bf g} + (\varepsilon - 1) \nabla \cdot {\bf v}_{\bf g} = 0$ %where $Q:= \{ q \in L^2_0(\Omega_{\varepsilon}), \; q=0 \mbox{ in } \Omega \setminus \bar{\Omega}_{\varepsilon} \}$ 
is supplemented by the relation $(\varepsilon - 1) \nabla \cdot {\bf v}_{\bf g} + p_{\bf g} = 0$ in $\Omega$. %In short, the convexity of $\Omega$ strongly supports \eqref{regularity}.\\
Akin to the classical Stokes system, the operator associated with this system is weakly coercive over $[H^1(\Omega) \cap L^2_0(\Omega)]^3 \times L^2_0(\Omega)$ equipped with the natural norm. This can be verified by choosing in the underlying variational form a test field ${\bf w}=(\varepsilon-1/2){\bf v}_{\bf g}$ in the first equation, and a test-function $q=p_{\bf g}$ in the second equation. Thus the convexity of $\Omega$ strongly supports \eqref{regularity}.\\
Now in order to establish the consistency of the explicit scheme \eqref{eq7} we first introduce an auxiliary field $\hat{\bf e}_h(\cdot,t)$  
belonging to ${\bf V}_h$ for every $t \in [0,T]$, uniquely defined up to an additive field depending only on $t$ as follows, for every $t \in [0,T]$:

\begin{equation} 
\label{consist1}
(\nabla \hat{\bf e}_h\{\cdot,t\}, \nabla {\bf v}) + (\nabla \cdot \hat{\bf e}_h\{\cdot,t\}, \nabla \cdot {\bf v})_{\varepsilon-1} =
(\nabla {\bf e}\{\cdot,t\}, \nabla {\bf v}) + (\nabla \cdot {\bf e}\{\cdot,t\}, \nabla \cdot {\bf v})_{\varepsilon-1} \; \forall {\bf v} \in {\bf V}_h. 
\end{equation}
The time-dependent additive field up to which $\hat{\bf e}_h\{\cdot,t\}$ is defined can be determined by requiring that 
$\int_{\Omega} \{\hat{\bf e}_h(\cdot,t)-{\bf e}(\cdot,t)\}d{\bf x} = {\bf 0}$ $\forall t \in [0,T]$.\\

Let us further assume that for every $t \in [0,T]$ ${\bf e}(\cdot,t) \in [H^2(\Omega)]^3$. In this case, from classical approximation results based 
on the interpolation error, we can assert that,
\begin{equation}
\label{consist1bis} \| \nabla \{ \hat{\bf e}_h(\cdot,t) - {\bf e}(\cdot,t)\} \| + \| \nabla \cdot \{ \hat{\bf e}_h(\cdot,t) - {\bf e}(\cdot,t)\} \|_{\varepsilon-1} \leq \hat{C}_1 h \| {\mathcal H}({\bf e})(\cdot,t) \| \, \forall t \in [0,T], 
\end{equation}
where $\hat{C}_1$ is a mesh-independent constant.\\
Let us show that there exists another mesh-independent constant $\hat{C}_0$ such that for every $t \in [0,T]$ it holds,
\begin{equation}
\label{consist0}
 \| \hat{\bf e}_h(\cdot,t) - {\bf e}(\cdot,t) \| \leq \hat{C}_0 h^2 \| {\mathcal H}({\bf e})(\cdot,t) \| \, \forall t \in [0,T].
\end{equation}
Since $\hat{\bf e}_h(\cdot,t) - {\bf e}(\cdot,t) \in [L^2_0(\Omega)]^3$ for every $t$, we may write:
\begin{equation}
\label{consist2}
\| \hat{\bf e}_h(\cdot,t) - {\bf e}(\cdot,t) \| = \displaystyle \sup_{{\bf g} \in  [L^2_0(\Omega)]^3 \setminus \{{\bf 0}\}} 
\frac{(\hat{\bf e}_h\{\cdot,t\} - {\bf e}\{\cdot,t\},{\bf g})}{\| {\bf g} \|} \; \forall t \in [0,T].
\end{equation}
Defining ${\bf W} := \{ {\bf w} \; | \; {\bf w} \in [H^2(\Omega) \cap L^2_0(\Omega)]^3 \; \partial_n {\bf w} = {\bf 0} \mbox{ on } \partial \Omega \}$, 
owing to \eqref{regularity} we have $\forall t \in [0,T]$:
\begin{equation}
\label{consist3}
\| \hat{\bf e}_h(\cdot,t) - {\bf e}(\cdot,t) \| \leq C \displaystyle \sup_{{\bf w} \in  {\bf W} \setminus \{{\bf 0}\}} 
\frac{(\hat{\bf e}_h\{\cdot,t\} - {\bf e}\{\cdot,t\}, -\nabla^2 {\bf w} - \nabla \{\varepsilon - 1\} \nabla \cdot {\bf w})}{\| {\mathcal H}({\bf w}) \|}.
\end{equation}
Since $\partial_n {\bf w} = {\bf 0}$ and $\varepsilon =1$ on $\partial \Omega$ we may integrate by parts the numerator in \eqref{consist3} to obtain for every 
$t \in [0,T]$,
\begin{equation}
\label{consist4}
\| \hat{\bf e}_h(\cdot,t) - {\bf e}(\cdot,t) \| \leq C \displaystyle \sup_{{\bf w} \in  {\bf W} \setminus \{{\bf 0}\}} 
\frac{(\nabla [\hat{\bf e}_h\{\cdot,t\} - {\bf e}\{\cdot,t\}], \nabla {\bf w}) + (\nabla \cdot [\hat{\bf e}_h\{\cdot,t\} - {\bf e}\{\cdot,t\}],
\nabla \cdot {\bf w})_{\varepsilon - 1}}{\| {\mathcal H}({\bf w}) \|}.
\end{equation}
Taking into account (\ref{consist1}) the numerator of \eqref{consist4} can be rewritten as 
\[ (\nabla[\hat{\bf e}_h\{\cdot,t\} - {\bf e}\{\cdot,t\}], \nabla [ {\bf w} - {\bf v}]) + (\nabla \cdot [\hat{\bf e}_h\{\cdot,t\} - {\bf e}\{\cdot,t\}],
\nabla \cdot [{\bf w}-{\bf v}])_{\varepsilon - 1} \; \forall {\bf v} \in {\bf V}_h. \]
Then choosing ${\bf v}$ to be the ${\bf V}_h$-interpolate of ${\bf w}$, taking into account \eqref{consist1bis} we easily establish \eqref{consist0} 
with $\hat{C}_0 = C \hat{C}_1^2$.\\
To conclude these preliminary considerations, we refer to Chapter 5 of \cite{Ruas}, to conclude that the second order time-derivative $\partial_{tt} \hat{\bf e}_h(\cdot,t)$ is well defined in $[H^1(\Omega)]^3$ for every $t \in [0,T]$, as long as $\partial_{tt} {\bf e}(\cdot,t)$ lies in $[H^1(\Omega)]^3$ for every $t \in [0,T]$. 
Moreover, provided $\partial_{tt} {\bf e} \in [H^2(\Omega)]^3$ for every $t \in [0,T]$, the following estimate holds:
\begin{equation}
\label{consist0bis}
 \| \partial_{tt}\{\hat{\bf e}_h(\cdot,t) - {\bf e}(\cdot,t)\} \| \leq \hat{C}_0 h^2 \| {\mathcal H}(\partial_{tt}{\bf e})(\cdot,t) \| \; \forall t \in [0,T].
\end{equation} 
In addition to the results given in this sub-section, we recall that, according to the Sobolev Embedding Theorem, there exists a constant 
$C_T$ depending only on $T$ such that it holds:
\begin{equation}
\label{Sobolemb} 
\| u \|_{L^{\infty}(0,T)} \leq C_T [ \| u \|_{L^2(0,T)}^2 + \| u_t \|_{L^2(0,T)}^2 ]^{1/2} \; \forall u \in H^1(0,T).
\end{equation}

In the remainder of this work we assume a certain regularity of ${\bf e}$, namely, \\

\noindent \underline{\textit{Assumption}$^*$ :} The solution ${\bf e}$ to equation \eqref{eq1} belongs to $[H^4\{\Omega \times(0,T)\}]^3$.  \\

\noindent Now taking $u = \| {\mathcal H}({\bf e})(\cdot,t) \|$ %= [\int_{\Omega} | {\mathcal H}({\bf e})({\bf x},t)|^2 d{\bf x}]^{1/2}$,  
we have $u_t(t) = \displaystyle \frac{\int_{\Omega} \{{\mathcal H}({\bf e}) : {\mathcal H}(\partial_t{\bf e})\}({\bf x},t) d{\bf x}}{\| \{{\mathcal H}({\bf e})\}(\cdot,t)\|}$, where $:$ denotes the inner product of two constant tensors of order greater than or equal to three. Then by the Cauchy-Schwarz inequality and taking into account \textit{Assumption} $^*$, it trivially follows from \eqref{Sobolemb} that the following upper bound holds:
\begin{equation}
\label{hessianbound}
 \| \{{\mathcal H}({\bf e})\}(\cdot,t) \| \leq C_T [ \| \| {\mathcal H}({\bf e}) \| \|_{L^2(0,T)}^2 + 
\| \| {\mathcal H}(\partial_t{\bf e}) \| \|_{L^2(0,T)}^2 ]^{1/2} \; \forall t \in (0,T).
\end{equation}

In complement to the above ingredients we extend the inner products $(\cdot,\cdot)_{\varepsilon_h,h}$ and $(\cdot,\cdot)_{\partial \Omega,h}$, and associated norms  
$\parallel \cdot \parallel_{\varepsilon_h,h}$ and $\parallel \cdot \parallel_{\partial \Omega,h}$ in a semi-definite manner, to fields in $[L^2(\Omega)]^N$, $N \leq 3$ as follows: \\
First of all, $\forall K \in {\mathcal T}_h$ let $\Pi_K : L^2(K) \rightarrow P_1(K)$ be the standard orthogonal projection operator onto the space $P_1(K)$ of linear functions in $K$. We set 

\[ \forall {\bf u} \mbox{ and }{\bf v} \in [L^2(\Omega)]^N, \; ({\bf u},{\bf v})_{\varepsilon_h,h} := \displaystyle \sum_{K \in {\mathcal T}_h} 
\frac{\varepsilon(G_K) volume(K)}{4} \sum_{i=1}^4 \Pi_K{\bf u}(S_{K,i}) \cdot \Pi_K{\bf v}(S_{K,i}). \]

Let us generically denote by $F \subset \partial \Omega$ a face of tetrahedron $K$ such that $area(K \cap \partial \Omega)>0$. Moreover we denote by $\Pi_F(v)$ the standard orthogonal projection of a function $v \in L^2(F)$ onto the space of linear functions on $F$. Similarly we define:
  
\[ \forall {\bf u} \mbox{ and }{\bf v} \in [L^2(\partial \Omega)]^N, \; ({\bf u},{\bf v})_{\partial \Omega,h} := \displaystyle \sum_{F \subset \partial \Omega} 
\frac{area(F)}{3} \sum_{i=1}^3 \Pi_F{\bf u}(R_{F,i}) \cdot \Pi_F{\bf v}(R_{F,i}). \]

It is noteworthy that whenever ${\bf u}$ and ${\bf v}$ belong to ${\bf V}_h$, both semi-definite inner products coincide with the inner products previously defined for 
such fields.\\
The following results hold in connection to the above inner products: 
\begin{lemma}
\label{JCAM} 
Let $\Delta_{\varepsilon_h}({\bf u},{\bf v}):=({\bf u},{\bf v})_{\varepsilon_h,h}-({\bf u},{\bf v})_{\varepsilon_h}$ for ${\bf u},{\bf v} \in [L^2(\Omega)]^N$    
There exists a mesh independent constant $c_{\Omega}$ such that $\forall {\bf u} \in [H^1(\Omega)]^3$ and $\forall {\bf v} \in {\bf V}_h$, 
\begin{equation}
\label{Delta}
|\Delta_{\varepsilon_h}({\bf u},{\bf v})| \leq c_{\Omega} \parallel \varepsilon \parallel_{0,\infty} h \parallel \nabla{\bf u} \parallel \parallel {\bf v} \parallel_h 
\mbox{ \rule{2mm}{2mm}}.
\end{equation}
\end{lemma} 
\begin{lemma}
\label{JCAMbis}
Let $ \Gamma_{h}(\gamma\{{\bf u}\},\gamma\{{\bf v}\}):=(\gamma\{{\bf u}\},\gamma\{{\bf v}\})_{\partial \Omega,h}-(\gamma\{{\bf u}\},\gamma\{{\bf v}\})_{\partial \Omega}$ for ${\bf u},{\bf v} \in [H^1(\Omega)]^N$, where $\gamma\{w\}$ represents the trace on $\partial \Omega$ of a function $w \in H^1(\Omega)$. Let also 
$\nabla_{\partial \Omega}$ be the tangential gradient operator over $\partial \Omega$.   
There exists a mesh independent constant $\tilde{c}_{\partial \Omega}$ such that $\forall {\bf u} \in [H^{2}(\Omega)]^3$ and $\forall {\bf v} \in {\bf V}_h$,
\begin{equation}
\label{tildeGamma} 
|\Gamma_{h}(\gamma\{{\bf u}\},\gamma\{{\bf v}\}) | \leq \tilde{c}_{\partial \Omega} h \parallel \nabla_{\partial \Omega} \gamma\{{\bf u}\} \parallel_{\partial \Omega}\parallel \gamma\{{\bf v}\} \parallel_{\partial \Omega,h}. \mbox{ \rule{2mm}{2mm}}.
\end{equation}
\end{lemma} 
The proof of Lemma \ref{JCAM} is based on the Bramble-Hilbert Lemma, and we refer to \cite{JCAM2010} for more details. Lemma \ref{JCAMbis} in turn 
follows from the same arguments combined with the Trace Theorem, which ensures that $\nabla_{\partial \Omega} \gamma\{{\bf u}\}$ is well defined in $[L^2(\partial \Omega)]^3$ if ${\bf u} \in [H^2(\Omega)]^3$. Incidentally the Trace Theorem allows us to bound above the right hand side of \eqref{tildeGamma} in such a way 
that the following estimate also holds for another mesh independent constant $c_{\partial \Omega}$:
\begin{equation}
\label{Gamma} 
|\Gamma_{h}(\gamma\{{\bf u}\},\gamma\{{\bf v}\}) | \leq c_{\partial \Omega} h 
[\parallel \nabla {\bf u} \parallel^2 + \parallel {\mathcal H}({\bf u}) \parallel^2]^{1/2} \parallel \gamma\{{\bf v}\} \parallel_{\partial \Omega,h}. \mbox{ \rule{2mm}{2mm}}.
\end{equation} 
To conclude we prove the validity of the following upper bounds:
\begin{lemma}
\label{JCAMter}  
$\forall {\bf v} \in [L^2(\Omega)]^3$ it holds  
$$\|{\bf v}\|_{\varepsilon_h,h} \leq \sqrt{5} \parallel {\bf v} \parallel_{\varepsilon_h}.$$ 
\end{lemma}
\underline{Proof:} Denoting by $\Pi_h({\bf v})$ the function defined in $\Omega$ whose restriction to every $K \in {\mathcal T}_h$ is $\Pi_K({\bf v})$ for a given  
${\bf v} \in [L^2(\Omega)]^3$, from an elementary property of the orthogonal projection we have 
\begin{equation} 
\label{intermezzo}
\parallel \Pi_{h}({\bf v}) \parallel_{\varepsilon_h} \leq \parallel {\bf v} \parallel_{\varepsilon_h}, \; \forall {\bf v} \in [L^2(\Omega)]^3. 
\end{equation}
Now taking $v$ such that $v_{|K} \in P_1(K) \; \forall K \in {\mathcal T}_h$, by a straightforward calculation using the expression of $v_{|K}$ in terms of barycentric coordinates we have:
\[ \int_K \varepsilon(G_K) v_{|K}^2 d{\bf x} = \displaystyle \frac {\varepsilon(G_K) volume(K)}{10} \left\{ \displaystyle \sum_{i=1}^4 |v(S_{K,i})|^2 + 
 \displaystyle \sum_{i=1}^3 \displaystyle \sum_{j=i+1}^4 v(S_{K,i})v(S_{K,j}) \right\}. \]
It trivially follows that 
\[ \int_K \varepsilon(G_K) v_{|K}^2 d{\bf x} = \displaystyle \frac {\varepsilon(G_K) volume(K)}{20} \left\{ \displaystyle \sum_{i=1}^4 |v(S_{K,i})|^2 + 
\displaystyle \left[ \sum_{i=1}^4 v(S_{K,i}) \right]^2 \right\}. \]
and finally 
\[ 5 \int_K \varepsilon(G_K) v_{|K}^2 d{\bf x} \geq \displaystyle \frac {\varepsilon(G_K) volume(K)}{4} \displaystyle \sum_{i=1}^4 |v(S_{K,i})|^2. \] 
This immediately yields Lemma \ref{JCAMter}, taking into account \eqref{intermezzo}. \rule{2mm}{2mm}
\begin{lemma}
\label{JCAMqua}  
$\forall {\bf v} \in [L^2(\partial \Omega)]^3$ it holds 
$\|{\bf v}\|_{\partial \Omega,h} \leq 2 \parallel {\bf v} \parallel_{\partial \Omega}$. 
\rule{2mm}{2mm}
\end{lemma}
%Noticing that the trace on $\partial \Omega$ of a function in $H^1(\Omega)$ is well-defined in $L^2(\partial \Omega)$, 
The proof of this Lemma is based on arguments  entirely analogous to Lemma \ref{JCAMter}.     
 
\subsection{Residual estimation}
To begin with we define for $k=0,1,\ldots,N$ functions $\hat{\bf e}^k_h \in {\bf V}_h$ by $\hat{\bf e}^k_h(\cdot) := \hat{\bf e}_h(\cdot,k \tau)$.
In the sequel for any function or field $A$ defined in $\Omega \times (0,T)$, $A^k(\cdot)$ denotes 
$A(\cdot,k \tau)$, except for other quantities carrying the subscript $h$ such as ${\bf e}_h^k$.\\

Let us substitute ${\bf e}^k_h$ by $\hat{\bf e}^k_h$ for 
$k=2,3,\ldots,N$ on the left hand side of the first equation of \eqref{eq7} and take also as initial conditions $\hat{\bf e}_h^j$ instead of ${\bf e}_h^j$, $j=0,1$.\\
The case of the initial conditions will be dealt with in the next section in the framework of the convergence analysis. As for   
the variational residual $E_h^k({\bf v})$ resulting from the above substitution, where $E_h^k$ being a linear functional acting on ${\bf V}_h$, 
it can be expressed in the following manner:
\begin{equation}
\label{consist7}
\begin{array}{l}
E_h^k({\bf v}) = ( \{\partial_{tt}{\bf e}\}^k,{\bf v})_{\varepsilon} + (\nabla {\bf e}^k,\nabla {\bf v}) + (\nabla \cdot {\bf e}^k,\nabla \cdot {\bf v})_{\varepsilon-1} + (\nabla \varepsilon \cdot {\bf e}^k,\nabla \cdot {\bf v}) + ( \{\partial_t {\bf e}\}^k, {\bf v})_{\partial \Omega} \\ 
\\
+ \bar{\bf F}^k({\bf v}) + (\bar{d}^k,\nabla \cdot {\bf v})+ \bar{\bf G}^k({\bf v}) \; \forall {\bf v} \in {\bf V}_h,
\end{array} 
\end{equation} 
%\displaystyle ( \varepsilon \{ {\bf T}^k_h({\bf e}) + {\bf R}^k({\bf e}) \} + 
% + ({\bf S}^k({\bf e})+{\bf Q}^k({\bf e}),{\bf v})_{\partial \Omega} \\
\noindent where  
\begin{equation}
\label{residual2}
\begin{array}{l}
\bar{d}^k = \nabla \varepsilon \cdot \{ \hat{\bf e}^k_h - {\bf e}^k \},
\end{array}
\end{equation}
\noindent and $\bar{\bf F}^k({\bf v})$ and $\bar{\bf G}^k({\bf v})$ can be written as follows:
\begin{equation}
\label{residual1}
\begin{array}{l}
\bar{\bf F}^k({\bf v})= {\bf F}^k_1({\bf v})+{\bf F}^k_2({\bf v})+{\bf F}^k_3({\bf v})+{\bf F}^k_4({\bf v}),\\
\mbox{with} \\
{\bf F}_1^k({\bf v})= ({\bf T}^k(\hat{\bf e}_h-{\bf e}),{\bf v})_{\varepsilon_h,h};\\
{\bf F}_2^k({\bf v})= \Delta_{\varepsilon_h} ({\bf T}^k({\bf e}),{\bf v} );\\
{\bf F}_3^k({\bf v})= ((\varepsilon_h-\varepsilon){\bf T}^k({\bf e}),{\bf v});\\
{\bf F}_4^k({\bf v})= ({\bf T}^k({\bf e}) - \{\partial_{tt}{\bf e}\}^k,{\bf v})_{\varepsilon}, 
%\\{\bf Q}^k_h({\bf e}) := \displaystyle \frac{\{\hat{\bf e}^{k+1}_h-{\bf e}^{k+1}\} - \{\hat{\bf e}^{k-1}_h-{\bf e}^{k-1}\}}{2 \tau}, \\
\end{array}
\end{equation}
\noindent ${\bf T}^k$ being the finite-difference operator defined by,
\begin{equation}
\label{Tk}
{\bf T}^k(\cdot) := \displaystyle \frac{\{\cdot\}^{k+1} - 2 \{\cdot\}^{k} +\{\cdot\}^{k-1}}{\tau^2},
\end{equation}
\noindent and
\begin{equation}
\label{residual3}
\begin{array}{l}
\bar{\bf G}^k({\bf v})= {\bf G}^k_1({\bf v})+{\bf G}^k_2({\bf v})+{\bf G}^k_3({\bf v}),\\
\mbox{with} \\
{\bf G}_1^k({\bf v})= ({\bf Q}^k(\hat{\bf e}_h-{\bf e}),{\bf v})_{\partial \Omega,h};\\
{\bf G}_2^k({\bf v})= \Gamma_{h}({\bf Q}^k({\bf e}),{\bf v}) ;\\
%+ \displaystyle \left(\frac{{\bf e}^{k+1} - 2 {\bf e}^{k} + {\bf e}^{k-1}}{\tau^2},{\bf v} \right)_{\varepsilon_h-\varepsilon} 
{\bf G}_3^k({\bf v})= ({\bf Q}^k({\bf e})-\{\partial_t{\bf e}\}^k,{\bf v} )_{\partial \Omega},
\end{array}
\end{equation}
\noindent ${\bf Q}^k$ being the finite-difference operator defined by,
\begin{equation}
\label{Qk}
{\bf Q}^k(\cdot) := \displaystyle \frac{\{\cdot\}^{k+1} - \{\cdot\}^{k-1}}{2 \tau}.
%\\{\bf Q}^k_h({\bf e}) := \displaystyle \frac{\{\hat{\bf e}^{k+1}_h-{\bf e}^{k+1}\} - \{\hat{\bf e}^{k-1}_h-{\bf e}^{k-1}\}}{2 \tau}, \\
\end{equation}

Notice that, under \textit{Assumption} $^*$ both $\partial_{t} {\bf e}(\cdot,t)$ and $\partial_{tt} {\bf e}(\cdot,t)$ belong to $[H^1(\Omega)]^3$ for every $t \in [0,T]$. Hence we can define $\partial_{t} \hat{\bf e}_h$ from $\partial_{t} {\bf e}$ and 
$\partial_{tt} \hat{\bf e}_h$ from $\partial_{tt} {\bf e}$ in the same way as $\hat{\bf e}_h$ is defined from ${\bf e}$. Moreover straightforward calculations lead to,
\begin{equation}
\label{consist80}
  {\bf T}^k({\cdot}) = \displaystyle \frac{1}{\tau^2} \left\{ \int_{(k-1) \tau}^{k \tau} \left( \int_t^{ k \tau} \partial_{tt} \{\cdot\} \; ds 
\right) dt	+ \displaystyle \int_{ k \tau}^{(k+1) \tau} \left( \int_{k \tau}^t \partial_{tt} \{\cdot\} \; ds \right)dt \right\}.
\end{equation}
Furthermore another straightforward calculation allows us writing:
\begin{equation}
\label{consist81}
\begin{array}{l}
  {\bf T}^k({\bf e}) = \{\partial_{tt}{\bf e}\}^k + {\bf R}^k({\bf e}) \\
	\\
	\mbox{where }\\
	{\bf R}^k({\bf e}) := \displaystyle \frac{1}{2 \tau^2} \left[ -\int_{(k-1) \tau}^{k \tau} \{ t - (k-1)\tau \}^2 \partial_{ttt} {\bf e} dt + \int_{k \tau}^{(k+1) \tau} \{ (k+1)\tau - t \}^2 \partial_{ttt} {\bf e} dt \right].
	\end{array}
	\end{equation}
Similarly,
\begin{equation}
\label{consist82}
  {\bf Q}^k(\cdot) = \displaystyle \frac{1}{2\tau} \int_{(k-1) \tau}^{(k+1) \tau} \partial_{t} \{\cdot\} \; dt,  
%+ \displaystyle \int_{ k \tau}^{(k+1) \tau} \left( \partial_{t} \hat{\bf e}_h - \partial_{t} {\bf e} \right)d \right\}.
\end{equation}
and
\begin{equation}
\label{consist83}
\begin{array}{l}
{\bf Q}^k({\bf e}) = \{\partial_{t}{\bf e}\}^k + {\bf S}^k({\bf e}), \\
\mbox{where } \\
{\bf S}^k({\bf e}) = \displaystyle \frac{1}{2 \tau} \left[ - \int_{(k-1) \tau}^{k \tau} \{ t - (k-1)\tau \} \partial_{tt} {\bf e} dt 
+ \int_{k \tau}^{(k+1) \tau} \{ (k+1)\tau - t \} \partial_{tt} {\bf e} dt \right].
\end{array}
\end{equation} 
		
Now we note that the sum of the terms on the first line of the expression of $E_h^k({\bf v})$ equals zero because they are just the left hand side of \eqref{eq2} 
at time $t = k \tau$. Therefore the functions $\hat{\bf e}_h^k \in {\bf V}_h$ are the solution of the following problem, for $k=1,2,\ldots,N-1$:
 
\begin{equation}\label{consist7bis}
  \begin{array}{l}
   \displaystyle  \left (\varepsilon \frac{\hat{\bf e}_h^{k+1} - 2 \hat{\bf e}_h^k + \hat{\bf e}_h^{k-1}}{\tau^2}, {\bf v} \right) +
    (\nabla \hat{\bf e}_h^k, \nabla {\bf v}) + (\nabla \cdot \{\varepsilon \hat{\bf e}_h^k\}, \nabla \cdot {\bf v}) 
    - (\nabla \cdot \hat{\bf e}_h^k, \nabla \cdot {\bf v}) \\
   + \displaystyle \left (\frac{\hat{\bf e}_h^{k+1} - \hat{\bf e}_h^{k-1}}{2\tau}, {\bf v} \right)_{\partial \Omega}= \bar{\bf F}^k({\bf v} ) + ( \bar{d}^k, \nabla \cdot {\bf v} ) + 
	\bar{\bf G}^k({\bf v}) ~ \forall {\bf v} \in {\bf V}_h,\\
    \\
    \hat{\bf e}_h^0(\cdot) = \hat{\bf e}_h(\cdot,0) \mbox{ and } \hat{\bf e}_h^1(\cdot) = \hat{\bf e}_h(\cdot,\tau) \mbox{ in } \Omega, 
  \end{array}
\end{equation}  
$\bar{d}^k$, $\bar{\bf F}^k$ and $\bar{\bf G}^k$ being given by \eqref{residual2}, \eqref{residual1}-\eqref{Tk} and \eqref{residual3}-\eqref{Qk}. \\
Estimating $\| \bar{d}^k \|$ is a trivial matter. Indeed, since ${\bf e}^k \in [H^2(\Omega)]^3$, from \eqref{consist0} we immediately obtain,
\begin{equation}
\label{Normdk}
\| \bar{d}^k \| \leq \hat{C}_0 h^2 | \varepsilon |_{1,\infty} \| {\mathcal H}({\bf e}^k) \|.
\end{equation}
\noindent Therefore consistency of the scheme will result from suitable estimations of $| \bar{\bf F}^k |_h$ and $| \bar{\bf G}^k |_{\partial \Omega,h}$ 
in terms of ${\bf e}$, $\varepsilon$, $\tau$ and $h$, which we next carry out.\\ 
%With this aim we assume that ${\bf e} \in [H^4\{\Omega \times (0,T)\}]^3$. \\
First of all we derive some upper bounds for the operators ${\bf T}^k(\cdot)$, ${\bf Q}^k(\cdot)$, ${\bf R}^k(\cdot)$ and ${\bf S}^k(\cdot)$. With this aim 
we denote by $| \cdot |$ the euclidean norm of $\Re^M$, for $M \geq 1$. \\
From \eqref{consist80} and the Cauchy-Schwarz inequality, we easily derive for every ${\bf x} \in \Omega$ and ${\bf u}$ 
such that $\{\partial_{tt}{\bf u}\}(\cdot,t) \in [H^{2}(\Omega)]^3$ $\forall t \in (0,T)$,
%From \eqref{consist81}, by straightforward upper bounds, followed by the Cauchy-Schwarz inequality, we obtain
\begin{equation*}
\begin{array}{l}
| {\bf T}^k_h\{{\bf u}({\bf x}) \} | \leq \displaystyle \frac{1}{\tau^2} \left\{ \int_{(k-1) \tau}^{k \tau} \left[ \int_t^{ k \tau} \left| \{\partial_{tt} {\bf u}\}({\bf x}) \right| ds 
\right]dt	+ \displaystyle \int_{ k \tau}^{(k+1) \tau} \left[ \int_{k \tau}^t \left| \{\partial_{tt} {\bf u}\}({\bf x}) \right|ds \right]dt \right\} \\
\leq \displaystyle \frac{1}{\tau^2} \left\{ \int_{(k-1) \tau}^{k \tau} \left[ \int_{(k-1)\tau}^{k \tau} \left|  \{\partial_{tt} {\bf u}\}({\bf x}) \right| ds 
\right]dt	+ \displaystyle \int_{ k \tau}^{(k+1) \tau} \left[ \int_{k \tau}^{(k+1)\tau} \left| \{ \partial_{tt} {\bf u}\}({\bf x}) \right|ds \right]dt \right\} \\
%= \displaystyle \frac{1}{\tau} \left[ \int_{(k-1)\tau}^{k \tau} \left| \partial_{tt} \hat{\bf u}_h - \partial_{tt} {\bf u} \right| dt	
%+ \displaystyle \int_{k \tau}^{(k+1)\tau} \left| \partial_{tt} \hat{\bf u}_h - \partial_{tt} {\bf u} \right|dt \right] =
= \displaystyle \frac{1}{\tau} \int_{(k-1)\tau}^{(k+1) \tau} \left| \{\partial_{tt} {\bf u}\}({\bf x}) \right| dt.   
\end{array}
\end{equation*}
It follows that, for every ${\bf u}$ such that $\{\partial_{tt}{\bf u}\}(\cdot,t) \in [H^{2}(\Omega)]^3$ $\forall t \in (0,T)$ we have,
\begin{equation}
\label{consist7qui}
| {\bf T}^k_h\{{\bf u}({\bf x}) \} | \leq \displaystyle \sqrt{\frac{2}{\tau}} \left[ \int_{(k-1)\tau}^{(k+1) \tau} \left| \{ \partial_{tt} {\bf u}\}({\bf x}) \right|^2 dt \right]^{1/2}.
\end{equation} 
Furthermore, from \eqref{consist81} and the inequality $a + b \leq [2(a^2+b^2)]^{1/2} \; \forall a,b \in \Re$, 
%and generically denoting ${\bf R}^k({\bf e}({\bf x})$ by ${\bf R}^k({\bf e})$ 
for every ${\bf x} \in \Omega$ we obtain:
\begin{equation*}
\begin{array}{l}
|{\bf R}^k\{{\bf e}({\bf x})\}| \leq \displaystyle \frac{1}{2} \left[ \int_{(k-1) \tau}^{k \tau} \displaystyle \left| \{\partial_{ttt} {\bf e}\}({\bf x}) \right| dt + \int_{k \tau}^{(k+1) \tau} \displaystyle \left| \{\partial_{ttt} {\bf e}\}({\bf x}) \right| dt \right] \\
\leq \displaystyle \frac{\sqrt{\tau}}{2} \left[ \left\{\int_{(k-1) \tau}^{k \tau} \displaystyle \left| \{\partial_{ttt} {\bf e}\}({\bf x}) \right|^2 dt \right\}^{1/2} + \left\{\int_{k \tau}^{(k+1) \tau} \displaystyle \left| \{\partial_{ttt} {\bf e}\}({\bf x}) \right|^2 dt \right\}^{1/2} \right].
\end{array}
\end{equation*}
It follows that,
\begin{equation}
\label{consist7ter}
|{\bf R}^k\{{\bf e}({\bf x})\}| \leq \displaystyle \sqrt{\frac{\tau}{2}} \left[ \int_{(k-1) \tau}^{(k+1) \tau} \displaystyle \left| \{\partial_{ttt} {\bf e}\}({\bf x}) \right|^2 dt \right]^{1/2}.
\end{equation}

On the other hand from \eqref{consist82} and the Cauchy-Schwarz inequality we trivially have for every ${\bf x} \in \partial \Omega$ and 
${\bf u}$ such that $\gamma\{\partial_{t}{\bf u}\}(\cdot,t) \in [H^{3/2}(\partial \Omega)]^3$ $\forall t \in (0,T)$:
\begin{equation}
\label{consist7sex}
| {\bf Q}^k\{{\bf u}({\bf x})\} | \leq \displaystyle \frac{1}{\sqrt{2\tau}} \left[\int_{(k-1) \tau}^{(k+1) \tau} |\{\partial_{t} {\bf u}\}({\bf x})|^2 dt 
\right]^{1/2}.
\end{equation}  
Finally, similarly to \eqref{consist7ter}, from \eqref{consist83} for every ${\bf x}$ in $\partial \Omega$ we successively derive,
\begin{equation*}
\begin{array}{l}
|{\bf S}^k\{{\bf e}({\bf x})\}| \leq \displaystyle \frac{1}{2} \left[ \int_{(k-1) \tau}^{k \tau} \displaystyle \left| \{\partial_{tt} {\bf e}\}({\bf x}) \right| dt + \int_{k \tau}^{(k+1) \tau} \displaystyle \left| \{\partial_{tt} {\bf e}\}({\bf x}) \right| dt \right] \\
\leq \displaystyle \frac{\sqrt{\tau}}{2} \left[ \left\{\int_{(k-1) \tau}^{k \tau} \displaystyle \left| \{\partial_{tt} {\bf e}\}({\bf x}) \right|^2 dt \right\}^{1/2} + \left\{\int_{k \tau}^{(k+1) \tau} \displaystyle \left| \{\partial_{tt} {\bf e}\}({\bf x}) \right|^2 dt \right\}^{1/2} \right].
\end{array}
\end{equation*}
Therefore it holds,
\begin{equation}
\label{consist7qua}
|{\bf S}^k\{{\bf e}({\bf x})\}| \leq \displaystyle \sqrt{\frac{\tau}{2}} \left[ \int_{(k-1) \tau}^{(k+1) \tau} \displaystyle \left| \{\partial_{tt} {\bf e}\}({\bf x}) \right|^2 dt \right]^{1/2}.
\end{equation}
Notice that bounds entirely analogous to \eqref{consist7qui} and \eqref{consist7sex} hold for $\nabla {\bf T}^k(\cdot)={\bf T}^k(\nabla \cdot)$ and 
$\nabla_{\partial \Omega} {\bf Q}^k(\cdot)={\bf Q}^k(\nabla_{\partial \Omega} \cdot)$, that is, $\forall {\bf x} \in \Omega$,
\begin{equation}
\label{consist7sep}
\begin{array}{l}
| \nabla {\bf T}^k_h({\bf e})({\bf x}) | \leq 
\displaystyle \sqrt{\frac{2}{\tau}} \left[ \int_{(k-1)\tau}^{(k+1) \tau} \left| \{ \partial_{tt} \nabla {\bf e}\}({\bf x}) \right|^2 dt \right]^{1/2},   
\end{array}
\end{equation}
\noindent and $\forall {\bf x} \in \partial \Omega$,
\begin{equation}
\label{consist7oct}
| \nabla_{\partial \Omega} {\bf Q}^k({\bf e})({\bf x}) | \leq \displaystyle \frac{1}{\sqrt{2\tau}} \left [\int_{(k-1) \tau}^{(k+1) \tau} |\{\partial_{t} 
\nabla_{\partial \Omega} {\bf e}\}({\bf x})|^2 dt \right]^{1/2}.
\end{equation} 

Next we estimate the four terms in the expression \eqref{residual1} of ${\bf F}^k({\bf v})$.\\
With the use of \eqref{consist7qui} and of Lemma \ref{JCAMter} followed by a trivial manipulation, we successively have: 
\begin{equation}
\label{F10}
\begin{array}{l}
|{\bf F}_1^k({\bf v})| \leq \parallel \varepsilon \parallel_{0,\infty} \displaystyle \sqrt{\frac{2}{\tau}} \left\|  \left[ \int_{(k-1)\tau}^{(k+1) \tau} \left| \{ \partial_{tt} (\hat{\bf e}_h-{\bf e})\}(\cdot,t) \right|^2 dt \right]^{1/2} \right\|_{h} \parallel {\bf v} \parallel_{h}\\
\leq \sqrt{5} \parallel \varepsilon \parallel_{0,\infty} \displaystyle \sqrt{\frac{2}{\tau}} \left\| \left[ \int_{(k-1)\tau}^{(k+1) \tau} \left| \{ \partial_{tt} (\hat{\bf e}_h-{\bf e})\}(\cdot,t) \right|^2 dt \right]^{1/2} \right\| \parallel {\bf v} \parallel_h \\
\leq \sqrt{5} \parallel \varepsilon \parallel_{0,\infty} \displaystyle \sqrt{\frac{2}{\tau}} \left[ \int_{(k-1)\tau}^{(k+1) \tau} \left\| \{ \partial_{tt} (\hat{\bf e}_h-{\bf e})\}(\cdot,t)  \right\|^2 dt \right]^{1/2} \parallel {\bf v} \parallel_h
\end{array}
\end{equation}
Recalling \eqref{residual1} and applying \eqref{consist0bis} to \eqref{F10}, we come up with,
\begin{equation}
\label{F1}
|{\bf F}_1^k({\bf v})| \leq \hat{C}_0 h^2 \displaystyle \sqrt{\frac{10}{\tau}} \parallel \varepsilon \parallel_{0,\infty} \left[ \int_{(k-1)\tau}^{(k+1) \tau} \parallel \{ {\mathcal H}(\partial_{tt}{\bf e})\}(\cdot,t) \parallel^2 dt \right]^{1/2} \parallel {\bf v} \parallel_h.
\end{equation}  
Next,  combining \eqref{Delta} and \eqref{consist7sep} we immediately obtain. 
\begin{equation}
\label{F2}
|{\bf F}_2^k({\bf v})| \leq c_{\Omega} h \parallel \varepsilon \parallel_{0,\infty} 
\displaystyle \sqrt{\frac{2}{\tau}} \left[ \int_{(k-1)\tau}^{(k+1) \tau} \left|\ \{ \partial_{tt} \nabla {\bf e}\}(\cdot,t) \right\|^2 dt \right]^{1/2} 
\parallel {\bf v} \parallel_h.
\end{equation} 
Further, from \eqref{upperbound}, \eqref{consist80}, \eqref{consist7sep} and the standard estimate $\parallel \varepsilon_h - \varepsilon \parallel_{0,\infty} \leq C_{\infty} h 
| \varepsilon |_{1,\infty}$ where $C_{\infty}$ is a mesh-independent constant, we derive  
\begin{equation}
\label{F3}
|{\bf F}_3^k({\bf v})| \leq C_{\infty} h 
| \varepsilon |_{1,\infty} \displaystyle \sqrt{\frac{2}{\tau}} \left[ \int_{(k-1)\tau}^{(k+1) \tau} \left\| \{ \partial_{tt} {\bf e}\}(\cdot,t) \right\|^2 dt \right]^{1/2} \parallel {\bf v} \parallel_h
\end{equation}  
Finally by \eqref{upperbound}, \eqref{consist81} and \eqref{consist7ter}, we have
\begin{equation}
\label{F4}
|{\bf F}_4^k({\bf v})| \leq \parallel \varepsilon \parallel_{0,\infty} \displaystyle \sqrt{\frac{\tau}{2}} \left[ \int_{(k-1) \tau}^{(k+1) \tau} \displaystyle \left\| \{\partial_{ttt} {\bf e}\}(\cdot,t) \right\|^2 dt \right]^{1/2} \parallel {\bf v} \parallel_h.
\end{equation}

Now we turn our attention to the three terms in the expression \eqref{residual3} of ${\bf G}^k({\bf v})$. 
First of all, owing to Assumption$^{*}$ and standard error estimates, we can write for a suitable mesh-independent constant $\hat{C}_1$:
\begin{equation}
\label{H1estimate} 
\parallel \{\nabla \partial_t(\hat{\bf e}_h - {\bf e})\}(\cdot,t) \parallel \leq \hat{C}_1 h \parallel \{{\mathcal H}(\partial_t{\bf e})\}(\cdot,t)) \parallel, \; \forall t \in [0,T].
\end{equation}
On the other hand, by the Trace Theorem there exists a contant $C_{Tr}$ depending only on $\Omega$ such that,
\begin{equation}
\label{trace} 
\parallel \{\partial_t(\hat{\bf e}_h - {\bf e})\}(\cdot,t) \parallel_{\partial \Omega} \leq C_{Tr}  [\parallel \{\partial_t (\hat{\bf e}_h - {\bf e})\}(\cdot,t) \parallel^2 + \parallel \{\nabla \partial_t (\hat{\bf e}_h - {\bf e})\}(\cdot,t) \parallel^2 ]^{1/2}\; \forall t \in [0,T].
\end{equation}
Hence by \eqref{consist0bis}, \eqref{H1estimate} and \eqref{trace}, we have, for a suitable mesh-independent constant $C_B$: 
\begin{equation}
\label{boundarestim}
\parallel \{\partial_t(\hat{\bf e}_h - {\bf e})\}(\cdot,t) \parallel_{\partial \Omega} \leq C_{B} h \parallel \{{\mathcal H}(\partial_t{\bf e})\}(\cdot,t) \parallel \; \forall t \in [0,T].
\end{equation}
Now recalling \eqref{residual1} and taking into account \eqref{consist7oct} and Lemma \ref{JCAMqua}, similarly to \eqref{F10} we first obtain:
\begin{equation}
\label{G10}
%\begin{array}{l}
| {\bf G}_1^k({\bf v}) | \leq \displaystyle \frac{2}{\sqrt{2\tau}} \left[ \int_{(k-1) \tau}^{(k+1) \tau} \left\| \{\partial_{t} (\hat{\bf e}_h-{\bf e})\}(\cdot,t) \right\|_{\partial \Omega}^2 dt \right]^{1/2} \parallel {\bf v} \parallel_{\partial \Omega,h}.
\end{equation}    
Then using \eqref{boundarestim} we immediately establish,
 \begin{equation}
\label{G1}
%\begin{array}{l}
| {\bf G}_1^k({\bf v}) | \leq C_B h  \displaystyle \sqrt{\frac{2}{\tau}} \left [\int_{(k-1) \tau}^{(k+1) \tau} \parallel \{{\mathcal H}(\partial_t{\bf e})\}(\cdot,t) \parallel^2 dt \right]^{1/2} \parallel {\bf v} \parallel_{\partial \Omega,h}.
\end{equation} 
Next we switch to ${\bf G}^k_2$. Using \eqref{Gamma} and \eqref{consist7oct}, similarly to \eqref{F2} we derive,
\begin{equation}
\label{G2}
|{\bf G}_2^k({\bf v})| \leq c_{\partial \Omega} h  
\displaystyle \sqrt{\frac{2}{\tau}} \left[ \int_{(k-1)\tau}^{(k+1) \tau} \left(\parallel \{\nabla \partial_t {\bf e}\}(\cdot,t) \parallel^2 + 
\parallel \{{\mathcal H}(\partial_t {\bf e})\}(\cdot,t) \parallel^2 \right) dt \right]^{1/2} \parallel {\bf v} \parallel_{\partial \Omega,h}.
\end{equation} 
As for $G_3^k$, taking into account \eqref{consist83} together with \eqref{consist7qua} and \eqref{upperboundary}, we obtain:
\begin{equation}
\label{G30}
|{\bf G}_3^k({\bf v})| \leq 
\displaystyle \sqrt{\frac{\tau}{2}} \left[ \int_{(k-1) \tau}^{(k+1) \tau} \displaystyle \left\| \{\partial_{tt} {\bf e}\}(\cdot,t) \right\|^2_{\partial \Omega} dt \right]^{1/2} 
\parallel {\bf v} \parallel_{\partial \Omega,h}.
\end{equation}
Then using the Trace Theorem (cf. \eqref{trace}), we finally establish,
\begin{equation}
\label{G3}
|{\bf G}_3^k({\bf v})| \leq 
 C_{Tr} \displaystyle \sqrt{\frac{\tau}{2}} \left[ \int_{(k-1) \tau}^{(k+1) \tau} \displaystyle \left( \parallel \{\partial_{tt} {\bf e}\}(\cdot,t) \parallel^2 + 
\parallel \{\nabla \partial_{tt} {\bf e}\}(\cdot,t) \parallel^2 \right) dt \right]^{1/2} 
\parallel {\bf v} \parallel_{\partial \Omega,h}.
\end{equation} 

Now collecting \eqref{F1}, \eqref{F2}, \eqref{F3} and \eqref{F4} we can write,
\begin{equation}
\label{NormFk}
\begin{array}{l}
| \bar{F}^k |_h \leq \hat{C}_0 h^2 \displaystyle \sqrt{\frac{10}{\tau}} \parallel \varepsilon \parallel_{0,\infty} \left[ \int_{(k-1)\tau}^{(k+1) \tau} \parallel \{ {\mathcal H}(\partial_{tt}{\bf e})\}(\cdot,t) \parallel^2 dt \right]^{1/2} \\ 
+ c_{\Omega} h \parallel \varepsilon \parallel_{0,\infty} 
\displaystyle \sqrt{\frac{2}{\tau}} \left[ \int_{(k-1)\tau}^{(k+1) \tau} \left|\ \{ \partial_{tt} \nabla {\bf e}\}(\cdot,t) \right\|^2 dt \right]^{1/2} \\
+ C_{\infty} h 
| \varepsilon |_{1,\infty} \displaystyle \sqrt{\frac{2}{\tau}} \left[ \int_{(k-1)\tau}^{(k+1) \tau} \left\| \{ \partial_{tt} \nabla {\bf e}\}(\cdot,t) \right\|^2 dt \right]^{1/2} \\
+ \parallel \varepsilon \parallel_{0,\infty} \displaystyle \sqrt{\frac{\tau}{2}} \left[ \int_{(k-1) \tau}^{(k+1) \tau} \displaystyle \left\| \{\partial_{ttt} {\bf e}\}(\cdot,t) \right\|^2 dt \right]^{1/2}.
\end{array}
\end{equation}
On the other hand \eqref{G1}, \eqref{G2} and \eqref{G3} yield,  
%Now collecting \eqref{F1}, \eqref{F2}, \eqref{F3} and \eqref{F4} we can write,
\begin{equation}
\label{NormGk}
\begin{array}{l}
| \bar{G}^k |_{\partial \Omega,h} \leq C_B h  \displaystyle \sqrt{\frac{2}{\tau}} \left [\int_{(k-1) \tau}^{(k+1) \tau} \parallel \{{\mathcal H}(\partial_t{\bf e})\}(\cdot,t) \parallel^2 dt \right]^{1/2} \\
+ c_{\partial \Omega} h  
\displaystyle \sqrt{\frac{2}{\tau}} \left[ \int_{(k-1)\tau}^{(k+1) \tau} \left(\parallel \{\nabla \partial_t {\bf e}\}(\cdot,t) \parallel^2 + 
\parallel \{{\mathcal H}(\partial_t {\bf e})\}(\cdot,t) \parallel^2 \right) dt \right]^{1/2} \\
+ C_{Tr} \displaystyle \sqrt{\frac{\tau}{2}} \left[ \int_{(k-1) \tau}^{(k+1) \tau} \displaystyle \left( \parallel \{\partial_{tt} {\bf e}\}(\cdot,t) \parallel^2 + 
\parallel \{\nabla \partial_{tt} {\bf e}\}(\cdot,t) \parallel^2 \right) dt \right]^{1/2}.
\end{array}
\end{equation}
Then, taking into account \eqref{consist7bis} and the stability condition 
\eqref{CFL}, by inspection we can assert that the consistency of scheme \eqref{eq7} is an immediate consequence of \eqref{Normdk}, \eqref{NormFk} and \eqref{NormGk}.
 
\section{Convergence results}

In order to establish the convergence of scheme \eqref{eq7} we combine the stability and consistence results obtained in the previous sections. With this aim we 
define $\bar{\bf e}^k_h := \hat{\bf e}^k_h- {\bf e}^k_h$ for $k=0,1,2,\ldots,N$. By linearity 
we can assert that the variational residual on the left hand side of the first equation of \eqref{eq7} for 
$k=2,3,\ldots,N$, when the ${\bf e}^k_h$s are replaced with the $\bar{\bf e}^k_h$s, and ${\bf e}_h^j$ is replaced with $\bar{\bf e}_h^j$ for $j=0,1$, is exactly 
$E_h^k({\bf v})$, since the residual corresponding to the ${\bf e}^k_h$'s vanishes by definition.
The initial conditions $\bar{\bf e}^j_h$ for $j=0$ and $j=1$ corresponding to the thus modified problem have to be estimated. This is the purpose of the next subsection.   

\subsection{Initial-condition deviations}
Here we turn our attention to the estimate of $\bar{E}_0$, which accounts for the deviation in the initial conditions appearing in the stability inequality \eqref{stability} that applies to the modification of \eqref{eq7} when ${\bf e}^k_h$ is replaced by $\bar{\bf e}^k_h$. \\
Let us first define,
\begin{equation}
\label{consist8bis}
  \begin{array}{l} 
	\hat{\bf e}_{1h} := (\hat{\bf e}^1_h - \hat{\bf e}^0_h)/\tau \\
	\\
	\bar{\bf e}_{1h} := \hat{\bf e}_{1h}-{\bf e}_{1h}, \\
	\\
	\bar{E}_0 := 
	\displaystyle \left\|\bar{\bf e}_{1h} \right\|^2_{\varepsilon_h,h} 
	+ \displaystyle \frac{\tau^2}{2} \left( 
  \|\nabla \bar{\bf e}_{1h} \|^2 + \| \nabla \cdot \bar{\bf e}_{1h} \|_{\varepsilon-1}^2 \right)  \\
		+ \displaystyle \frac{1}{2} \left( \| \nabla \bar{\bf e}_h^1 \|^2 +  \|\nabla  \bar{\bf e}_h^{0} \|^2 
	\right)
	+ \displaystyle \frac{1}{2} \left( \| \nabla \cdot \bar{\bf e}_h^1\|_{\varepsilon-1}^2  + \| \nabla \cdot \bar{\bf e}_h^{0} \|_{\varepsilon-1}^2 \right) 
	+ T | \varepsilon |_{2,\infty} \| \bar{\bf e}_h^0\|^2 
	\end{array}
\end{equation}
Recalling that ${\bf e}^1_h = {\bf e}^0_h + \tau {\bf e}_{1h}$ we have $\bar{\bf e}^1_h = \bar{\bf e}^0_h + \tau \bar{\bf e}_{1h}$. Thus, taking either 
${\mathcal A}=\nabla \bar{\bf e}_h^0$ or ${\mathcal A}=\nabla \cdot \bar{\bf e}_h^0$ and either ${\mathcal B} = \tau \nabla \bar{\bf e}_{1h}$ or 
${\mathcal B} = \tau \nabla \cdot \bar{\bf e}_{1h}$, 
we apply twice the inequality $\|{\mathcal A}+{\mathcal B}\|^2/2 \leq \| {\mathcal A} \|^2 + \|{\mathcal B}\|^2$ to \eqref{consist8bis} together with Lemma \ref{JCAMter}, to obtain, 
\begin{equation}
\label{consist8ter}
\begin{array}{l}
 \bar{E}_0 \leq  %\displaystyle \frac{\tau^2}{2} \left( \| \nabla \bar{\bf e}_{1h}\|^2  + \|\nabla \cdot \bar{\bf e}_{1h} \|^2 \right) 
\displaystyle \frac{3}{2} \left\{\|\nabla \bar{\bf e}_{h}^0 \|^2 + \| \nabla \cdot \bar{\bf e}_h^0\|_{\varepsilon-1}^2  
+ \tau^2 \left(\| \nabla \bar{\bf e}_{1h} \|^2 + \| \nabla \bar{\bf e}_{1h}\|^2_{\varepsilon-1} \right) \right\} \\
+ T | \varepsilon |_{2,\infty} \| \bar{\bf e}_h^0 \|^2 +  5 \| \varepsilon \|_{0,\infty} \| \bar{\bf e}_{1h} \|^2.
\end{array}
\end{equation}
Finally using the inequality $\| \nabla \cdot {\bf v} \|^2_{\varepsilon -1} \leq 3 \| \varepsilon - 1 \|_{0,\infty} \| \nabla {\bf v} \|^2$ $\forall {\bf v} \in 
[H^1(\Omega)]^3$, after straightforward manipulations we easily derive from \eqref{consist8ter}:   
\begin{equation}
\label{consist9bis}
\begin{array}{l}
 \bar{E}_0 \leq %\displaystyle \frac{\tau^2}{2} \left( \| \nabla \bar{\bf e}_{1h}\|^2  + \|\nabla \cdot \bar{\bf e}_{1h} \|^2 \right) 
\displaystyle \frac{3(1+ 3 \| \varepsilon - 1 \|_{0,\infty})}{2} \left( \|\nabla \bar{\bf e}_{h}^0 \|^2 + \tau^2 \| \nabla \bar{\bf e}_{1h} \|^2 \right) 
+ T | \varepsilon |_{2,\infty} \| \bar{\bf e}_h^0 \|^2 +  5 \| \varepsilon \|_{0,\infty} \| \bar{\bf e}_{1h} \|^2.
\end{array}
\end{equation}
We next use the obvious splitting $\bar{\bf e}_{h}^0 = (\hat{\bf e}_{h}^0 -{\bf e}_0) + ({\bf e}_0 -{\bf e}_{h}^0 )$, together with the trivial one, 
\begin{equation}
\label{consist10}
	\bar{\bf e}_{1h} = \displaystyle \frac{\hat{\bf e}_{h}^1 - \hat{\bf e}_{h}^0}{\tau}  - {\bf e}_{1h} = 
	 (\{\partial_t \hat{\bf e}_{h}\}_{|t=0} - {\bf e}_1 ) 
	+ \displaystyle \frac{1}{\tau} \int_0^{\tau} (\tau-t)\partial_{tt} \hat{\bf e}_h \; dt 
		+ ({\bf e}_1 - {\bf e}_{1h}). 
	  %+ \displaystyle \frac{1}{\tau} \int_0^{\tau} (\tau-t) \partial_{tt} {\bf e} \; dt.
 \end{equation}
Then plugging \eqref{consist10} into \eqref{consist9bis}, since $\left[ \sum_{i=1}^p \| A_i \| \right]^2 \leq p \sum_{i=1}^p \| A_i \|^2$  
for any set of $p$ functions or fields $A_i$, we obtain:
\begin{equation}
\label{consist11}
\begin{array}{l}
 \bar{E}_0 \leq \displaystyle \frac{3(1+ 3 \| \varepsilon - 1 \|_{0,\infty})}{2} \left[ 2 \left\{ \| \nabla (\hat{\bf e}_{h}^0 - {\bf e}_{0}) \|^2 + \| \nabla ({\bf e}_{0} - {\bf e}_{h}^0 ) \|^2 \right\} \right. \\
\left. 3 \tau^2 \displaystyle \left\{ \| \nabla \partial_t(\hat{\bf e}_h - {\bf e})(\cdot,0) \|^2 + \displaystyle \frac{1}{\tau^2} \left \|\int_0^{\tau} (\tau-t)\partial_{tt} \nabla \hat{\bf e}_h \; dt \right\|^2 + \| \nabla ({\bf e}_1 - {\bf e}_{1h}) \|^2 \right\} \right] \\ 
+ 2 T | \varepsilon \textit{}|_{2,\infty} \left( \| \hat{\bf e}_h^0 - {\bf e}_0 \|^2 +  \| {\bf e}_0 - {\bf e}_h^0 \|^2 \right) \\ 
+ 15 \| \varepsilon \|_{0,\infty} \displaystyle \left\{ \| \partial_t (\hat{\bf e}_h - {\bf e})(\cdot,0) \|^2 + \displaystyle \frac{1}{\tau^2} \left \|\int_0^{\tau} (\tau-t)\partial_{tt} \hat{\bf e}_h \; dt \right\|^2 + \| {\bf e}_1 - {\bf e}_{1h} \|^2 \right\}.
\end{array}
\end{equation}
Owing to a trivial upper bound and to the Cauchy-Schwarz inequality we easily derive
\begin{equation*} 
\begin{array}{l}
\displaystyle \frac{1}{\tau^2} \left\| \int_0^{\tau} (\tau-t)  \nabla \partial_{tt} \hat{\bf e}_h \; dt \right\|^2 \leq 
\left\| \int_0^{\tau} |\nabla \partial_{tt} \hat{\bf e}_h | dt \right\|^2 
\leq \tau \left\| \displaystyle \left[ \int_0^{\tau} | \nabla \partial_{tt} \hat{\bf e}_h |^2 dt \right]^{1/2} \right\|^2.
\end{array}
\end{equation*}
It follows that,
\begin{equation}
\label{consist12}
\displaystyle \frac{1}{\tau^2} \left\| \int_0^{\tau} (\tau-t)  \nabla \partial_{tt} \hat{\bf e}_h \; dt \right\|^2 
\leq  
\tau \int_0^{\tau} \left\| \nabla \partial_{tt} \hat{\bf e}_h \right\|^2 dt \leq (1+3 \| \varepsilon -1 \|_{0,\infty} ) \tau 
\int_0^{\tau} \left\| \nabla \partial_{tt} {\bf e} \right\|^2 dt.
\end{equation} 
The last inequality in \eqref{consist12} follows from the definition of $\hat{\bf e}_h$. Indeed we know that, 
\begin{equation}
\label{consist13}
(\nabla \partial_{tt} \hat{\bf e}_h,\nabla {\bf v}) + (\nabla \cdot \partial_{tt} \hat{\bf e}_h,\nabla \cdot {\bf v})_{\varepsilon-1} =
(\nabla \partial_{tt} {\bf e},\nabla {\bf v}) + (\nabla \cdot \partial_{tt} {\bf e},\nabla \cdot {\bf v})_{\varepsilon-1} \; \forall {\bf v} \in {\bf V}_h.
\end{equation}
Taking ${\bf v} = \hat{\bf e}_h$, by the Cauchy-Schwarz inequality, we easily obtain:
\begin{equation}
\label{consist14}
[ \| \nabla \partial_{tt} \hat{\bf e}_h \|^2  + \| \nabla \cdot \partial_{tt} \hat{\bf e}_h \|_{\varepsilon-1}^2 ]^{1/2} \leq 
[ \| \nabla \partial_{tt} {\bf e} \|^2  + \| \nabla \cdot \partial_{tt} {\bf e} \|_{\varepsilon-1}^2 ]^{1/2},
\end{equation} 
or yet 
\begin{equation}
\label{consist15}
\| \nabla \partial_{tt} \hat{\bf e}_h \| \leq \sqrt{1 + 3 \| \varepsilon -1 \|_{0,\infty}} \| \nabla \partial_{tt} {\bf e} \|,
\end{equation}
which validates \eqref{consist12}.\\
On the other hand according to \eqref{consist0bis} we have $\| \partial_{tt} \hat{\bf e}_h \| \leq \| \partial_{tt}{\bf e} \| + \hat{C}_0 
h^2 \| {\mathcal H}(\partial_{tt}{\bf e}) \|$. This yields,
\begin{equation}
\label{consist12bis}
\left\| \int_0^{\tau} (\tau-t) \partial_{tt} \hat{\bf e}_h \; dt \right\|^2 
\leq  
\tau \int_0^{\tau} \left\| \partial_{tt} \hat{\bf e}_h \right\|^2 dt \leq 2 \tau \left[\int_0^{\tau} \left\| \partial_{tt} {\bf e} \right\|^2 dt + \hat{C}_0^2 h^4  
\int_0^{\tau} \left\| {\mathcal H}(\partial_{tt} {\bf e}) \right\|^2 dt \right].
\end{equation} 
Incidentally we point out that \textit{Assumption} $^*$ and \eqref{Sobolemb} allow us to assert that
\begin{itemize}
\item
$\partial_{tt} \nabla {\bf e} \in [L^2\{\Omega \times (0,T)\}]^3$;
\item
$\| \partial_{tt} {\bf e} \| \in [ L^{\infty}(0,T)]^3$;
\item
$\{\partial_t {\bf e}\}_{|t=0} = {\bf e}_1 \in [H^2(\Omega)]^3$;
\item
$\{{\bf e}\}_{|t=0} = {\bf e}_0 \in [H^2(\Omega)]^3$.
\end{itemize}
Clearly enough, besides \eqref{consist0} and \eqref{consist1bis}, we will apply to \eqref{consist11} standard estimates based on the interpolation error in Sobolev norms (cf. \cite{Ciarlet}), together with the following obvious variants of \eqref{consist1bis} and \eqref{consist0bis}, namely,
\begin{equation}
\label{consist0ter} 
\begin{array}{l} 
\left\| \nabla \partial_{t} (\hat{\bf e}_h - {\bf e})(\cdot,0) \right\| \leq \hat{C}_1 h \left\| {\mathcal H}({\bf e}_1) \right\|, \\ 
\left\| \partial_{t} (\hat{\bf e}_h - {\bf e})(\cdot,0) \right\| \leq \hat{C}_0 h^2 \left\| {\mathcal H}({\bf e}_1) \right\|. 
\end{array}
\end{equation} 
Then taking into account that $\tau \leq 1/(2 \eta)$, from \eqref{consist11}-\eqref{consist12}-\eqref{consist12bis}-\eqref{consist0ter} 
and Assumption$^*$, we conclude that there exists a constant $\bar{C}_0$ depending on $\Omega$, $T$ and $\varepsilon$, but neither on $h$ nor on $\tau$, such that,
\begin{equation}
\label{E00}
\begin{array}{l}
\bar{E}_0 \leq \bar{C}_0^2 \displaystyle \left[(h^2 + \tau h^2 + h^4)\| {\mathcal H}({\bf e}_0) \|^2 + (h^2 + \tau^2 h^2)\| {\mathcal H}({\bf e}_1) \|^2 
+ \tau^3  \| \| \nabla \partial_{tt} {\bf e} \| \|_{L^2(0,T)}^2 \right. \\
\left. + \tau^2 \| \| \partial_{tt} {\bf e} \| \|_{L^{\infty}(0,T)}^2 + \tau h^4 \int_0^{\tau} \left\| {\mathcal H}(\partial_{tt} {\bf e}) \right\|^2 dt \right].
\end{array}
\end{equation}
Notice that, starting from \eqref{Sobolemb} with $u = \| \nabla \partial_{tt} {\bf e} \|$, similarly to \eqref{hessianbound}, we obtain
\[  
 \| \| \partial_{tt}{\bf e} \| \|_{L^{\infty}(0,T)} \leq C_T [ \| \| \partial_{tt}{\bf e} \| \|_{L^2(0,T)}^2 + 
\| \| \partial_{ttt}{\bf e}) \| \|_{L^2(0,T)}^2 ]^{1/2}.
\] 
Thus noting that $h \leq diam(\Omega)$, using again the upper bound $\tau \leq 1/(2 \eta)$ and extending the integral to the whole interval $(0,T)$ in \eqref{E00},  
from the latter inequality we infer the existence of another constant $C_0$ independent of $h$ and $\tau$, such that, 
\begin{equation}
\label{E0}
\begin{array}{l}
\bar{E}_0 \leq C_0^2 \displaystyle \left[h^2 \{\| {\mathcal H}({\bf e}_0) \|^2 + \| {\mathcal H}({\bf e}_1) \|^2 + \| {\mathcal H}(\partial_{tt} {\bf e}) \|^2_{L^2[\Omega \times (0,T)]} \} \right. \\
\left. + \tau^2 \{ \| \partial_{tt} {\bf e} \|^2_{L^2[\Omega \times (0,T)]} + \| \partial_{ttt} {\bf e} \|^2_{L^2[\Omega \times (0,T)]} 
+ \| \nabla \partial_{tt} {\bf e} \|^2_{L^2[\Omega \times (0,T)]} \} \right].
\end{array}
\end{equation}

\subsection{Error estimates}
In order to fully exploit the stability inequality \eqref{stability} we further define,
\begin{equation}
\label{BarSN}
	\bar{S}_N := \displaystyle \sum_{k=1}^{N-1}  \left( \displaystyle \frac{\tau}{2} | \bar{\bf F}^k |_h^2 + \displaystyle \frac{3}{2 \tau} \| \bar{d}^k\|^2 + 
	\displaystyle \frac{\tau}{2} \left| \bar{\bf G}^k \right|^2_{\partial \Omega,h} \right). 
\end{equation}
According to \eqref{consist7bis}, in order to estimate $\bar{S}_N$ under the regularity Assumption$^*$ on ${\bf e}$, we resort 
to the estimates \eqref{Normdk}, \eqref{NormFk} and \eqref{NormGk}. Using the inequality $|{\bf a} \cdot {\bf b}| \leq |{\bf a}| |{\bf b}|$ for ${\bf a},{\bf b} \in \Re^M$, it is easy to see that there exists two constants $\bar{C}_F$ and $\bar{C}_G$ independent of $h$ and $\tau$ such that
\begin{equation}
\label{BarFk}
\displaystyle \sum_{k=1}^{N-1} \displaystyle \frac{\tau}{2} | \bar{\bf F}^k |_h^2 \leq \bar{C}_F \displaystyle \left( h^2 + \tau^2 \right) | {\bf e} |_{H^4[\Omega \times (0,T)]}^2,
\end{equation}
\begin{equation}
\label{BarGk}
	\displaystyle \sum_{k=1}^{N-1} \displaystyle \frac{\tau}{2} | \bar{\bf G}^k |_{\partial \Omega,h}^2 \leq \bar{C}_G \displaystyle \left( h^2 + \tau^2 \right) | {\bf e} |_{H^3[\Omega \times (0,T)]}^2.
\end{equation}
%${\bf R}^k_h({\bf e})$, ${\bf T}^k_h({\bf e})$ and ${\bf S}^k_h({\bf e})$ given by \eqref{consist5}, \eqref{consist8} and \eqref{consist6}, respectively. \\
On the other hand, recalling \eqref{hessianbound} we have $k=1,2,\ldots,N-1$: 
\begin{equation}
\label{Hek}
\| {\mathcal H}({\bf e}^k) \|^2 \leq C_T^2 \displaystyle \int_0^T \left[ \| {\mathcal H}({\bf e}) \|^2 + \| {\mathcal H}(\partial_t{\bf e}) \|^2 \right] dt.
\end{equation}
Therefore, since $N = T/\tau$ we have:
\begin{equation}
\label{SHek}
\displaystyle \sum_{k=1}^{N-1} \| {\mathcal H}({\bf e}^k) \|^2 \leq \displaystyle \frac{C_T^2 T}{\tau} \displaystyle \int_0^T \left[ \| {\mathcal H}({\bf e}) \|^2 + \| {\mathcal H}(\partial_t{\bf e}) \|^2 \right] dt
\end{equation}
It follows from \eqref{SHek} and \eqref{Normdk} that,
\begin{equation}
\label{Bardk}
	\displaystyle \sum_{k=1}^{N-1} \displaystyle \frac{3}{2\tau} \| \bar{\bf d}^k \|^2 \leq \hat{C}_0^2 C_T ^2 T | \varepsilon |_{1,\infty}^2 \displaystyle \frac{3h^4}{2\tau^2} \left( | {\bf e} |_{H^2[\Omega \times (0,T)]}^2 + | {\bf e} |_{H^3[\Omega \times (0,T)]}^2 \right)
\end{equation}
Plugging \eqref{BarFk}, \eqref{Bardk} and \eqref{BarGk} into \eqref{BarSN} we can assert that there exists a constant $\tilde{C}$ depending on $\Omega$, $T$ and  
$\varepsilon$, but neither on $h$ nor on $\tau$, such that,
\begin{equation}
\label{SN}
\begin{array}{l}
\bar{S}_N \leq \tilde{C}^2 \displaystyle \left( h^2 + \tau^2 + \frac{h^4}{\tau^2} \right) \| {\bf e} \|_{H^4[\Omega \times (0,T)]}^2 
\end{array}
\end{equation}

\noindent Now recalling \eqref{stability}-\eqref{eq26}, provided \eqref{CFL} holds, together with $\tau \leq 1/(2 \eta)$, we have: 
\begin{equation}
\label{convergence0}
 \displaystyle   \left[ \left\| \frac{\bar{\bf e}_h^{m+1} - \bar{\bf e}_h^m}{\tau} \right\|^2 + \| \nabla \bar{\bf e}_h^{m+1}\|^2 + \| \nabla \bar{\bf e}_h^{m}\|^2 
\right]^{1/2} \leq 2 \sqrt{\bar{S}_N + 
		\bar{E}_0} e^{\beta/2}.
\end{equation}
This implies that, for $m=1,2,\ldots,N-1$, it holds:
\begin{equation}
\label{convergence1}
\begin{array}{l}
   \displaystyle  \left[ \left\| \frac{{\bf e}_h^{m+1} - {\bf e}_h^{m}}{\tau} - 
	\frac{{\bf e}^{m+1} - {\bf e}^{m}}{\tau} \right\|^2 + \| \nabla \left({\bf e}_h^{m+1}-{\bf e}^{m+1} \right) \|^2 
	+ \| \nabla \left({\bf e}_h^{m} - {\bf e}^m) \right) \|^2 \right]^{1/2} \\
	\leq  \displaystyle  \left[ \left\| \frac{\hat{\bf e}_h^{m+1} - {\bf e}^{m+1}}{\tau} - 
	\frac{\hat{\bf e}_h^{m} - {\bf e}^{m}}{\tau} \right\|^2 + \| \nabla \left(\hat{\bf e}_h^{m+1}-{\bf e}^{m+1} \right) \|^2 
	+ \| \nabla \left(\hat{\bf e}_h^{m} - {\bf e}^m \right) \|^2 \right]^{1/2} \\
	+\; 2 \sqrt{\bar{S}_N +	\bar{E}_0} e^{\beta/2}.
	\end{array}
\end{equation} 
Let us define a function ${\bf e}_h$ in $\bar{\Omega} \times [0,T]$ whose value at $t=k\tau$ equals ${\bf e}_h^k$ for $k=1,2,\ldots,N$ and that varies linearly with $t$ in each time interval $([k-1] \tau,k \tau)$, in such a way that $\partial_t{\bf e}_{h}({\bf x},t)=\displaystyle \frac{{\bf e}^k_h({\bf x}) - {\bf e}^{k-1}_h({\bf x})}{\tau}$ for every ${\bf x} \in \bar{\Omega}$ and $t \in ([k-1]\tau,k \tau)$. \\
Now we define $\tilde{A}^{m-1/2}(\cdot)$ for any function or field $A(\cdot,t)$ to be the mean value of $A(\cdot,t)$ in $(m \tau, [m+1]\tau$, that is 
$\tilde{A}^{m+1/2} = \tau^{-1} \int_{m \tau}^{(m+1)\tau} A(\cdot,t) dt$. Clearly enough we have
%\begin{equation}
%\label{convergence}\frac{{\bf e}^{m+1} - {\bf e}^m}{\tau} = \int_{m \tau}^{(m+1)\tau} {\bf e}_t dt
\begin{equation} 
\label{convergence2}
\displaystyle \left\| \frac{{\bf e}_h^{m+1} - {\bf e}_h^{m}}{\tau} - \frac{{\bf e}^{m+1} - {\bf e}^{m}}{\tau} \right\|^2 = 
 \displaystyle  \left\| \{\widetilde{\partial_t {\bf e}_h}\}^{m+1/2} - \{\widetilde{\partial_t {\bf e}}\}^{m+1/2} \right\|^2. 
\end{equation}
and also recalling \eqref{SHek} and \eqref{hessianbound} 
\begin{equation*}
\begin{array}{l}
\displaystyle \left\| \frac{\hat{\bf e}_h^{m+1} - {\bf e}^{m+1}}{\tau} - \frac{\hat{\bf e}_h^{m} - {\bf e}^{m}}{\tau} \right\|^2 
= \displaystyle \left\| \displaystyle \int_{m \tau}^{(m+1)\tau} \frac{\partial_t \hat{\bf e}_h - \partial_t {\bf e}}{\tau} dt \right\|^2 \\
\leq \left\| \displaystyle \frac{ \left[\int_{m \tau}^{(m+1)\tau} (\partial_t \hat{\bf e}_h - \partial_t {\bf e})^2 dt \right]^{1/2}}{\sqrt{\tau}}  \right\|^2 
= \displaystyle \frac{1}{\tau} \int_{m \tau}^{(m+1)\tau} \left\| \partial_t \hat{\bf e}_h - \partial_t {\bf e} \right\|^2 dt \\
\leq \displaystyle \frac{\hat{C}_0^2 h^4}{\tau} \int_{m \tau}^{(m+1)\tau} \| {\mathcal H}(\partial_t {\bf e}) \|^2 dt \leq 
C_T^2 \hat{C}_0^2 h^4 \int_0^T [\| {\mathcal H}(\partial_t {\bf e}) \|^2 + \| {\mathcal H}(\partial_{tt} {\bf e}) \|^2 ] dt,  
\end{array}
\end{equation*}
which implies that
\begin{equation}
\label{convergence3}
\begin{array}{l}
\displaystyle \left\| \frac{\hat{\bf e}_h^{m+1} - {\bf e}^{m+1}}{\tau} - \frac{\hat{\bf e}_h^{m} - {\bf e}^{m}}{\tau} \right\|^2 \leq C_T^2 \hat{C}_0^2 h^2
|diam(\Omega)|^2 \| {\bf e} \|_{H^4[ \Omega \times (0,T)]}^2.
\end{array}
\end{equation}
On the other hand from \eqref{consist1bis} and \eqref{hessianbound} we have for $j=m$ or $j=m+1$:
\begin{equation*}
\| \nabla (\hat{\bf e}_h^{j}-{\bf e}^{j} ) \|^2 \leq \hat{C}_1^2 h^2 \| {\mathcal H}({\bf e}^j) \|^2 \leq \hat{C}_1^2 C_T^2 h^2 \displaystyle 
\int_0^T \left[ \| {\mathcal H}({\bf e}) \|^2 +  \| {\mathcal H}(\partial_t {\bf e}) \|^2 \right] dt.
\end{equation*}
which yields,
\begin{equation}
\label{convergence4}
\| \nabla (\hat{\bf e}_h^{j}-{\bf e}^{j}) \|^2 \leq C_T^2 \hat{C}_1^2 h^2 \| {\bf e} \|_{H^3[ \Omega \times (0,T)]}^2.
\end{equation}
Now using Taylor expansions about $t=(m+1/2)\tau$ together with some arguments already exploited in this work, it is easy to establish that 
for a certain constant 
$\tilde{\mathcal C}$ independent of $h$ and $\tau$ it holds, 
\begin{equation}
\label{brokenpartialt}
\left\| \{ \widetilde{\partial_t {\bf e}}\}^{m+1/2} - [\partial_t {\bf e}]^{m+1/2} \right\| \leq \tilde{\mathcal C} \tau^2 
\left\| {\bf e} \right\|_{H^4[\Omega \times (0,T)]}, 
\end{equation}
where for every function $g$ defined in $\bar{\Omega} \times [0,T]$ and $s \in \Re^{+}$, $g^s$ is the function defined in $\Omega$ 
by $g^s(\cdot) = g(\cdot,s \tau)$.\\ 
Noticing that $\{\widetilde{\partial_t {\bf e}_h}\}^{m+1/2}(\cdot)$ is nothing but $[\partial_t {\bf e}_h]^{m+1/2}(\cdot)$, collecting \eqref{convergence1}, \eqref{convergence2}, \eqref{convergence3}, \eqref{convergence4}, \eqref{brokenpartialt}, 
together with \eqref{E0} and \eqref{SN}, we have thus proved the following convergence result for scheme \eqref{eq7}:\\

\textbf{Provided the CFL condition \eqref{CFL} is satisfied and $\tau$ also satisfies $\tau \leq 1/[2\eta]$, 
under Assumption $^*$ on ${\bf e}$, there exists a constant ${\mathcal C}$ depending only on $\Omega$, $\varepsilon$ and $T$ such that}
  
\begin{equation}
\label{convergence}
\boxed{
\begin{array}{l}
\displaystyle \max_{1 \leq m \leq N-1}  \left\| [\partial_t ({\bf e}_h - {\bf e})]^{m+1/2} \right\| 
+ \displaystyle \max_{2 \leq m \leq N}  \| \nabla ({\bf e}_h^{m} - {\bf e}^{m}) \| \\
\leq {\mathcal C} (\tau + h + h^2/\tau) \displaystyle \left\{ \| {\bf e} \|_{H^4[\Omega \times (0,T)]} + \| {\mathcal H}({\bf e}_0) \| + \| {\mathcal H}({\bf e}_1) \| \right\}. 
\mbox{ \rule{2mm}{2mm}}
\end{array}
}
\end{equation}   
\textbf{In short, as long as $\tau$ varies linearly with $h$, first order convergence of scheme \eqref{eq7} in terms of either $\tau$ or $h$ is thus established in the sense of the norms on the left hand side of \eqref{convergence}}. 

\section{Assessment of the scheme}
 
We performed numerical tests for the model problem \eqref{eq1}, taking $T=0.5$ and $\Omega$ to be the unit disk given by 
$\Omega= \{(x_1;x_2) |\;  r < 1\}$, where $r:= \sqrt{x_1^2+x_2^2}$.

We consider that the exact solution ${\bf e}$ of \eqref{eq1} is given by:

\begin{center}
${\bf e}=(e_1,e_2)$ with $e_1=-x_2 v_{\theta}(r,t)$ and $e_2 = x_1 v_{\theta}(r,t)$,\\
\vspace{2mm}
where $\displaystyle v_{\theta}(r,t) = \frac{e^{r-2t}}{\varepsilon(r)}$ \\
\end{center} 
with 
\begin{equation}
\label{eps}
\varepsilon(r) = 1+(1-4r^2)^m(1-H_{1/2}(r)), \mbox{for an integer } m > 1, 
\end{equation}
$H_{1/2}$ being the Heaviside function.\\ 

The initial data ${\bf e}_0$ and ${\bf e}_1$ are given by

\begin{align*}
{\bf e}_0 &:= {\bf e}_{|t=0} = \displaystyle \frac{(-x_2,x_1)e^r}{\varepsilon(r)}, \\
{\bf e}_1 &:= [\partial_t{\bf e}]_{|t=0} = \displaystyle \frac{2 (x_2,-x_1) e^r}{\varepsilon(r)}.
  \end{align*}
The right hand side ${\bf f}$ in turn is given by,

\begin{align*}
f_{1} &:= \varepsilon \partial_{tt} e_1 - \Delta e_1 + \partial_{x_1} (\nabla \cdot {\bf e}) = -4 x_2 e^{r-2t} - \Delta e_1, \\
f_{2} &:= \varepsilon \partial_{tt} e_2 - \Delta e_2 + \partial_{x_2} (\nabla \cdot {\bf e}) = 4 x_1 e^{r-2t} - \Delta e_2, 
\end{align*}
with 
\begin{align*}
\Delta e_1 &= - \displaystyle \frac{3x_2}{r} v_{\theta}^{'} - x_2 v_{\theta}^{''}, \\
\Delta e_2 &=  \displaystyle \frac{3x_1}{r} v_{\theta}^{'} + x_1 v_{\theta}^{''}, \\ v_{\theta}^{'} &= \displaystyle \frac{\varepsilon - \varepsilon^{'}}{\varepsilon^2} e^{r-2t}, \\
v_{\theta}^{''} &= \displaystyle 
\frac{\varepsilon^2-2\varepsilon \varepsilon^{'}-\varepsilon \varepsilon^{''} + 2 (\varepsilon^{'})^2}{\varepsilon^3} e^{r-2t}, \\
\varepsilon^{'} &= -8mr(1-4r^2)^{m-1}(1-H_{1/2}(r)), \\
\varepsilon^{''} &= 8m(8m^2r^2-4r^2-1)(1-4r^2)^{m-2}(1-H_{1/2}(r)). \\ 
\end{align*}

Notice that ${\bf e}$ satisfies absorbing conditions on the boundary of the unit disk. Indeed, 
\begin{align*}
\partial_r e_1 = \displaystyle \frac{x_1}{r} \partial_{x_1} e_1 + \displaystyle \frac{x_2}{r} \partial_{x_2} e_1 
= -\displaystyle \frac{x_2}{r} v_{\theta} - x_2 v_{\theta}^{'}
\end{align*}
and
\begin{align*}
\partial_r e_2 = \displaystyle \frac{x_{1}}{r} \partial_{x_1} e_2 + \displaystyle \frac{x_2}{r} \partial_{x_2} e_2 
= \displaystyle \frac{x_1}{r} v_{\theta} + x_1 v_{\theta}^{'}.
\end{align*}

Therefore, at $r=1$ we have
\begin{align*}
  \partial_n {\bf e} = \partial_r {\bf e}=(-x_2,x_1) (v_{\theta}(1)+v_{\theta}^{'}(1))=2(-x_2,x_1)e^{1-2t}.
  \end{align*}
Since ${\bf e} = (-x_2,x_1) e^{r-2t}$ for $r > 1/2$, it follows that $\partial_n {\bf e} + \partial_t {\bf e} = {\bf 0}$ for $r=1$.\\

In our computations we used the software package Waves \cite{waves} only
for the finite element method applied to the solution of the model
problem \eqref{eq1}. We note that this package was also used in \cite{cejm} to solve the
the same model problem \eqref{eq1} by a domain
decomposition FEM/FDM method.\\

We discretized the spatial domain $\Omega$ using a family of quasi-uniform meshes consisting of $2 \times 2^{2l+2}$ triangles $K$ for 
$l=1,...,6$, constructed as a certain mapping of the same number of triangles of the uniform mesh of the square $(-1,1)\times (-1,1)$, which is symmetric with respect to the cartesian axes and have their edges parallel to those axes and either to the line $x_1=x_2$ if $x_1x_2 \geq 0$ or to the line $x_1=-x_2$ otherwise. The above mapping is defined by means of a suitable transformation of cartesian into polar coordinates. For each value of $l$ we define a reference mesh size $h_l$ equal to $2^{-l}$.\\ 
We consider a partition of the time domain $(0,T)$ into time
 intervals $J=(t_{k-1},t_k]$ of equal length $\tau_l$
for a given number of time intervals $N, l=1,...,6$. 
We performed numerical tests taking $m=2,3,4,5$ in \eqref{eps}. We choose the time step
$\tau_l = 0.025 \times 2^{-l}, \; l=1,...,6$, which provides numerical stability for all meshes. 
%  such that
% $\tau = t_k - t_{k-1}$.
% As usual,
% we
% assume also a minimal angle condition of elements $K$ in the $K_h$.
We computed the maximum value over the time steps of the relative errors measured in the $L_2$-norms of the function, its gradient and its time-derivative in the polygons $\Omega_h$ defined to be the union of the triangles in the different meshes in use. Now ${\bf e}$ being the exact solution of \eqref{eq1}, ${\bf e}_h$ being the computed solution and setting $N := T /\tau_l$, these quantities are represented by,
\begin{equation}\label{norms}
  \begin{split}
e_l^1 &= \displaystyle \frac{\displaystyle \max_{1 \leq k \leq N} \| {\bf e}^k-{\bf e}_h^k \|_h}{\displaystyle \max_{1 \leq k \leq N}\| {\bf e}^k \|_h},\\
e_l^2 &= \displaystyle \frac{ \displaystyle \max_{1 \leq k \leq N} \| \nabla ({\bf e}^k - {\bf e}^k_h) \|_h}{
\displaystyle \max_{1 \leq k \leq N} \| \nabla {\bf e}^k \|_h}, \\
e_l^3 &=  \displaystyle \frac{\displaystyle \max_{1 \leq k \leq N-1} \| \{ \partial_t ({\bf e}-{\bf e}_h)\}^{k+1/2} \|_h}{\displaystyle \max_{1 \leq k \leq N-1} \| \{ \partial_t {\bf e}\}^{k+1/2}\|_h},
\end{split}
 \end{equation}
respectively, where $\| \cdot \|_h$ stands for the standard norm of $[L^2(\Omega_h)]^2$. Notice that, in the case of $P_1$ elements and of  convex domains, the error estimates in mean-square norms that hold for a polygonal domain extend to a curved domain, 
as long as the norm $\| \cdot \|$ is replaced by $\| \cdot \|_h$. \\   

In Tables 1-4 method's convergence in these three senses is observed taking $m=2,3,4,5$ in \eqref{eps}.
\begin{table}
\center
\begin{tabular}{ l  l l l l l l l l }
\hline
$l$ &  $nel$  & $nno$ &  $e_l^{(1)}$  & $e_{l-1}^{(1)}/e_{l}^{(1)}$ &
$e_l^{(2)}$
  & $e_{l-1}^{(2)}/e_{l}^{(2)}$ & $e_l^{(3)}$ &   $e_{l-1}^{(3)}/e_{l}^{(3)}$\\
\hline 
1 & 32    & 25     & 0.6057  &              & 2.2827  &           & 3.1375 &   \\
2 & 128   & 81     & 0.1499  &  4.0418      & 1.0769  &  2.1198   & 1.1536 & 2.7196 \\
3 & 512   & 289    & 0.0333  &  4.5007      & 0.4454  &  2.4178   & 0.5776 & 1.9972 \\
4 & 2048  & 1089   & 0.0078  &  4.2466      & 0.2077  &  2.1449   & 0.2802 & 2.0617 \\
5 & 8192  & 4225   & 0.0019  &  4.1288      & 0.1066  &  1.9483   & 0.1379 & 2.0313 \\
6 & 32768 & 16641  & 0.0005  &  4.0653      & 0.0535  &  1.9905  &  0.0690 &  1.9981 \\
\hline
\end{tabular}
\caption{ The computed maximum relative errors $e_n$ in maximum energy, maximum $L_2$ and  maximum in
  broken time error  on different  meshes with mesh sizes $h_l= 2^{-l}, l=1,...,6$ for the function $\varepsilon$
   with $m=2$ in  \eqref{eps}.}
\label{test1}
\end{table}

\begin{table}
\center
\begin{tabular}{ l  l l l l l l l l }
\hline
$l$ &  $nel$  & $nno$ &  $e_l^{(1)}$  & $e_{l-1}^{(1)}/e_{l}^{(1)}$ &
$e_l^{(2)}$
  & $e_{l-1}^{(2)}/e_{l}^{(2)}$ & $e_l^{(3)}$ &   $e_{l-1}^{(3)}/e_{l}^{(3)}$\\
\hline 
1 & 32    & 25     & 0.6144  &              & 1.8851 &           & 3.0462 &   \\
2 & 128   & 81     & 0.1511  &  4.0666      & 1.0794  &  1.7464   & 1.1417 & 2.6682 \\
3 & 512   & 289    & 0.0339  &  4.4553      & 0.4713  &  2.2904   & 0.5680 & 2.0099 \\
4 & 2048  & 1089   & 0.0080  &  4.2216      & 0.2166  &  2.1753   & 0.2760 & 2.0583 \\
5 & 8192  & 4225   & 0.0019  &  4.1207      & 0.1137  &  1.9049   & 0.1354 & 2.0381 \\
6 & 32768 & 16641  & 0.0005  &  4.0615      & 0.0566  &  2.0092  &  0.0677 &  1.9997 \\
\hline
\end{tabular}
\caption{ The computed maximum relative errors $e_n$ in maximum energy, maximum $L_2$ and  maximum in broken time error  on different  meshes with mesh sizes $h_l= 2^{-l}, l=1,...,6$ for the function $\varepsilon$
   with $m=3$ in  \eqref{eps}.}
\label{test2}
\end{table}

\begin{table}
\center
\begin{tabular}{ l  l l l l l l l l }
\hline
$l$ &  $nel$  & $nno$ &  $e_l^{(1)}$  & $e_{l-1}^{(1)}/e_{l}^{(1)}$ &
$e_l^{(2)}$
  & $e_{l-1}^{(2)}/e_{l}^{(2)}$ & $e_l^{(3)}$ &   $e_{l-1}^{(3)}/e_{l}^{(3)}$\\
\hline 
1 & 32    & 25     & 0.6122  &              & 1.9517 &           & 3.0545 &   \\
2 & 128   & 81     & 0.1529  &  4.0027      & 1.0896  &  1.7912   & 1.1445 & 2.6689 \\
3 & 512   & 289    & 0.0346  &  4.4266      & 0.4879  &  2.2331   & 0.5639 & 2.0296 \\
4 & 2048  & 1089   & 0.0082  &  4.2069      & 0.2234  &  2.1839   & 0.2728 & 2.0667 \\
5 & 8192  & 4225   & 0.0020  &  4.1151      & 0.1183  &  1.8879   & 0.1336 & 2.0418 \\
6 & 32768 & 16641  & 0.0005  &  4.0585      & 0.0595  &  1.9890  &  0.0668 &  2.0008 \\
\hline
\end{tabular}
\caption{ The computed maximum relative errors $e_n$ in maximum energy, maximum $L_2$ and  maximum in broken time error  on different  meshes with mesh sizes $h_l= 2^{-l}, l=1,...,6$ for the function $\varepsilon$
   with $m=4$ in  \eqref{eps}.}
\label{test3}
\end{table}

\begin{table}
\center
\begin{tabular}{ l  l l l l l l l l }
\hline
$l$ &  $nel$  & $nno$ &  $e_l^{(1)}$  & $e_{l-1}^{(1)}/e_{l}^{(1)}$ &
$e_l^{(2)}$
  & $e_{l-1}^{(2)}/e_{l}^{(2)}$ & $e_l^{(3)}$ &   $e_{l-1}^{(3)}/e_{l}^{(3)}$\\
\hline 
1 & 32    & 25     & 0.6107  &              & 1.9930 &            & 3.0603 &   \\
2 & 128   & 81     & 0.1546  &  3.9505      & 1.1006  &  1.8108   & 1.1464 & 2.6696 \\
3 & 512   & 289    & 0.0351  &  4.4031      & 0.4982  &  2.2090   & 0.5619 & 2.0403 \\
4 & 2048  & 1089   & 0.0084  &  4.1954      & 0.2288  &  2.1777   & 0.2706 & 2.0765 \\
5 & 8192  & 4225   & 0.0020  &  4.1106      & 0.1223  &  1.8715   & 0.1325 & 2.0417 \\
6 & 32768 & 16641  & 0.0005  &  4.0561     & 0.0607 &  2.0139     &  0.0662 &  2.0011 \\
\hline
\end{tabular}
\caption{ The computed maximum relative errors $e_n$ in maximum energy, maximum $L_2$ and  maximum in broken time error  on different  meshes with mesh sizes $h_l= 2^{-l}, l=1,...,6$ for the function $\varepsilon$
   with $m=5$ in  \eqref{eps}.}
\label{test4}
\end{table}

\begin{figure}[h!]
\begin{center}
\begin{tabular}{cc}
  {\includegraphics[ trim = 4.0cm 4.0cm 1.0cm 4.0cm, scale=0.5, clip=]{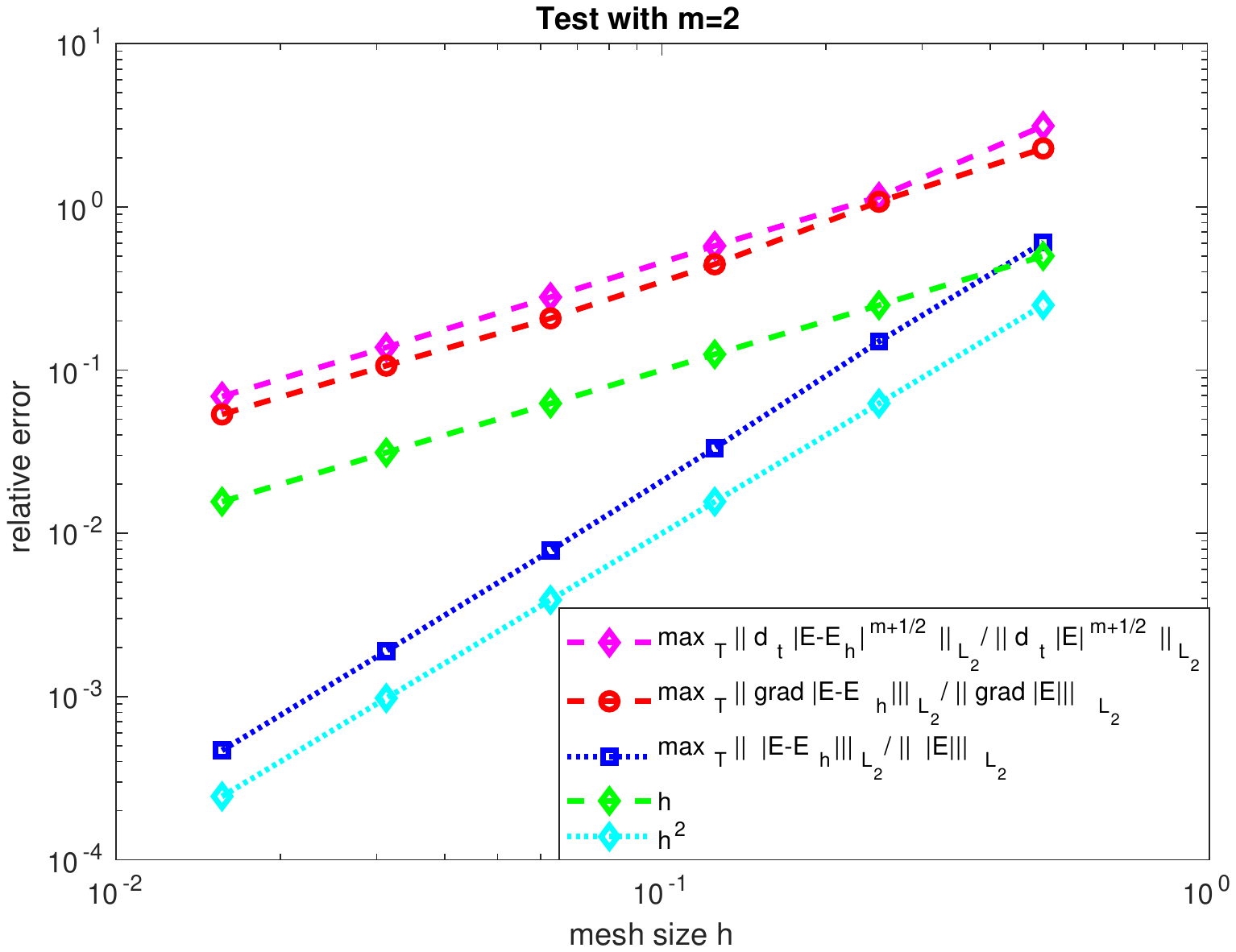}}  &
  {\includegraphics[trim = 4.0cm 4.0cm 1.0cm 4.0cm,   scale=0.5, clip=]{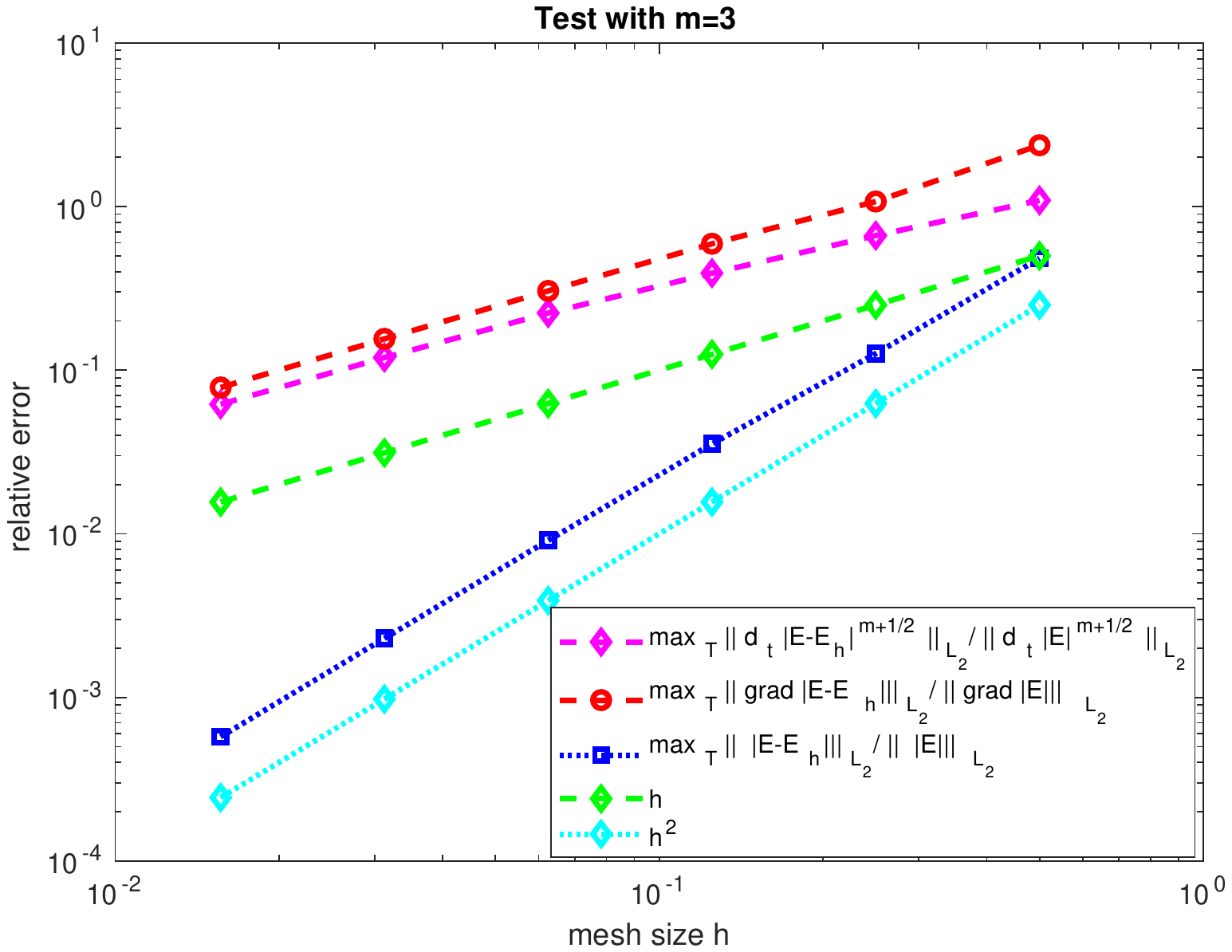}} \\
  $m=2$  & $m=3$ 
\end{tabular}
\end{center}
\caption{Maximum in time of relative errors for $m=2$ (left) and $m=3$ (right)}
\label{fig:F2}
\end{figure}

\begin{figure}[h!]
\begin{center}
\begin{tabular}{cc}
  {\includegraphics[trim = 4.0cm 4.0cm 1.0cm 4.0cm,  scale=0.5, clip=]{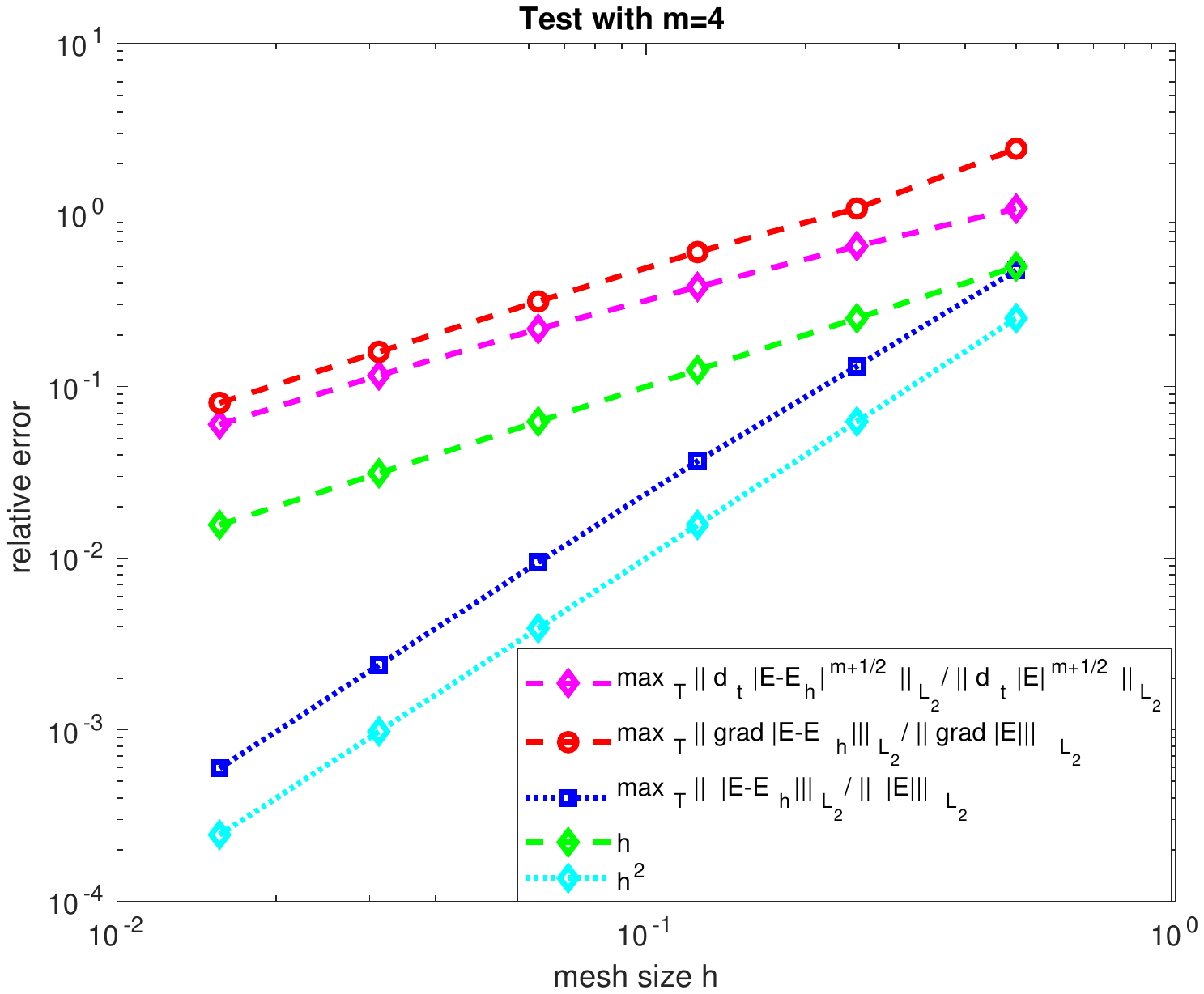}}  &
  {\includegraphics[trim = 4.0cm 4.0cm 1.0cm 4.0cm,  scale=0.5, clip=]{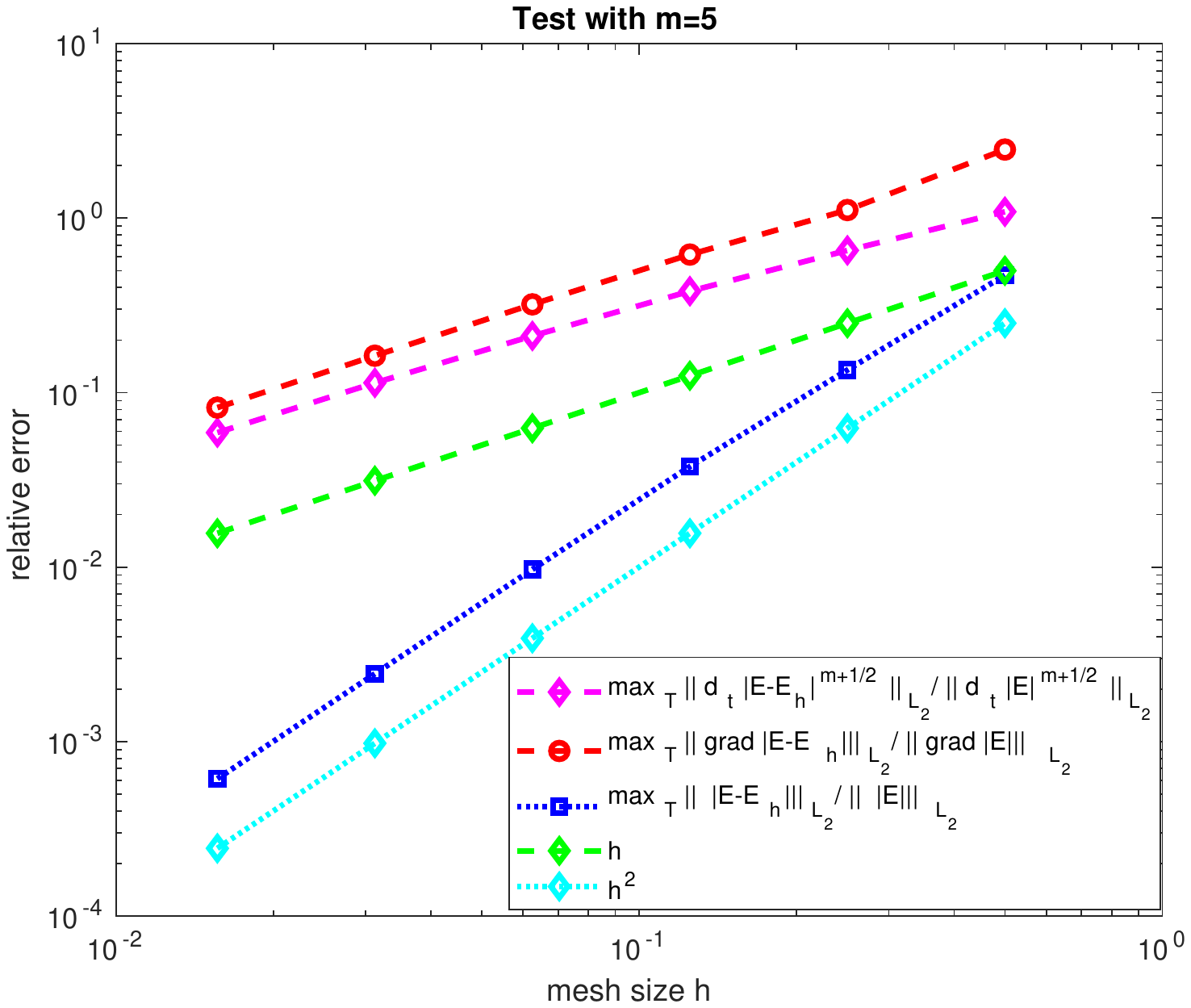}} \\
  $m=4$ & $m=5$ 
\end{tabular}
\end{center}
\caption{Maximum in time of relative errors for $m=4$ (left) and $m=5$ (right)}
\label{fig:F3}
\end{figure}

Figure \ref{fig:F2} shows convergence rates of our numerical scheme based on a $P_1$ space discretization, taking the 
function $\varepsilon$ defined by \eqref{eps} with $m=2$ (on the
left) and $m=3$ (on the right) for
$\varepsilon(r)$. Similar convergence results are presented in
Figures \ref{fig:F3} taking $m=4$ (on the left) and $m=5$ (on the right) in \eqref{eps}.

Observation of these tables and figures clearly indicates that our scheme behaves like a first order method in the (semi-)norm of $L^{\infty}[(0,T);H^1(\Omega)]$ for ${\bf e}$ and in the norm of $L^{\infty}[(0,T);L^2(\Omega)]$ for $\partial_t {\bf e}$ for all the chosen values of $m$.  
As far as the values of $m$ greater or equal to $4$ are concerned this perfectly conforms to the a priori error estimates 
established in Section 6. However those tables and figures also show that such theoretical predictions  
extend to cases not considered in our analysis such as $m=2$ and $m=3$, in which the regularity of the exact solution is lower than assumed. Otherwise stated some of our assuptions seem to be of academic interest only and a lower regularity of the solution such as $H^2[\Omega \times (0,T)]$ should be sufficient to attain optimal first order convergence in both senses.
On the other hand second-order convergence can be expected from our scheme in the norm of $L^{\infty}[(0,T);L^2(\Omega)]$ for ${\bf e}$, according to Tables 1-4 and Figures \ref{fig:F2} and \ref{fig:F3}. 

\section{Final remarks}

As previously noted, the approach advocated in this work was extensively and successfully tested in the framework of the solution of CIPs 
governed by Maxwell's equations. More specifically it was used with minor modifications to solve both the direct problem and the adjoint problem, as steps of an adaptive algorithm to determine the unknown electric permittivity. More details on this procedure can be found in \cite{MalmbergBeilina1,Malmberg}.\\

As a matter of fact the method studied in this paper was designed to handle composite dielectrics structured in such a way that layers with higher permittivity 
are completely surrounded by layers with a (constant) lower permittivity, say with unit value. It should be noted however that the assumption that the minimum value of $\varepsilon$ equal one in the outer layer was made here only to simplify things. Actually under the same assumptions \eqref{convergence} also applies to the case where $\varepsilon$ in inner layers is allowed to be smaller than in the outer layer, say $\varepsilon <1$. For instance, if $\varepsilon > 2/3$ the upper bound  
\eqref{convergence} also holds for a certain mesh-independent constant ${\mathcal C}$. This is because under such an assumption on $\varepsilon$ it is possible to guarantee that the auxiliary problems \eqref{consist1} and \eqref{vectpoisson} are coercive. \\
On the other hand in case $\varepsilon$ can be less than or equal to $2/3$, the convergence analysis of scheme \eqref{eq7} is a little more laborious. The key to the 
problem is a modification of the variational form \eqref{eq2} as follows. First of all we set $\varepsilon_{min} = \displaystyle \min_{{\bf x} \in \Omega} \varepsilon({\bf x})$. Then we recast \eqref{eq2} for every $t \in (0, T)$ as : $\forall  {\bf v} \in [H^1(\Omega)]^3$ it holds,  
\begin{equation}
\label{eq2bis}
%\begin{array}{ll}
 \left ( \varepsilon \partial_{tt} {\bf e},{\bf v} \right ) + (\nabla {\bf e},\nabla {\bf v})  +  \displaystyle \left(\nabla \cdot  \left[\frac{\varepsilon {\bf e}}
{\varepsilon_{min}}\right], \nabla \cdot {\bf v} \right) - 
(\nabla \cdot {\bf e}, \nabla \cdot {\bf v}) + (\partial_t{\bf e}, {\bf v})_{\partial \Omega}  = 0. %& (7.82)
%\end{array}
\end{equation}
Akin to \eqref{eq2}, problem (\ref{eq2bis}) is equivalent to Maxwell's equations \eqref{eq1}. Indeed,  
integrating by parts in \eqref{eq2bis}, for all $ {\bf v} \in [H^1(\Omega)]^3$ we get,
\begin{equation}\label{eq2ter}
\begin{array}{l}
\left ( \varepsilon \partial_{tt} {\bf e},{\bf v} \right ) + ( \nabla \times \nabla \times {\bf e}, {\bf v}) - (\nabla \nabla \cdot [\varepsilon {\bf e}],{\bf v}) \\
+ \left(\partial_n{\bf e} + \partial_t {\bf e},{\bf v} \right)_{\partial \Omega} +  (\nabla \cdot[\varepsilon_{min}^{-1} \varepsilon {\bf e}] 
- \nabla \cdot {\bf e}, {\bf v} \cdot {\bf n} )_{\partial \Omega}  = 0
\end{array}
\end{equation}
Besides the Maxwell's equations in $\Omega \times (0, T)$, this time the conditions on $\partial \Omega \times (0,T)$ are those resulting from \eqref{eq2ter}, 
that is,
\begin{equation}
\label{mixedbc}
\begin{array}{l}
(\partial_{n} {\bf e} + \partial_t {\bf e} ) \cdot {\bf n} + \nabla \cdot [\varepsilon_{min}^{-1} \varepsilon {\bf e}] 
- \nabla \cdot {\bf e} = 0, \\
(\partial_{n} {\bf e} + \partial_t {\bf e} ) \times {\bf n} = {\bf 0}.
\end{array}
\end{equation}
Since $\nabla \cdot \tilde{\bf e} = \nabla \cdot [\varepsilon \tilde{\bf e}] = 0$ on $\partial \Omega$, where $\tilde{\bf e}$ is 
the solution of Maxwell's equations in strong form \eqref{eq1}, this field necessarily satisfies the boundary conditions \eqref{mixedbc} as well.  
Then in the same way as the solution of \eqref{eq2}, this implies that the solution of \eqref{eq2bis} also fulfills $\nabla \cdot (\varepsilon {\bf e})=0$.\\
However in this case, given ${\bf g} \in [L^2_0(\Omega)]^3$, the auxiliary problem \eqref{vectpoisson} must be modified into,
%, instead of zero normal-derivative boundary conditions, the stationary counterpart of the boundary conditions \eqref{mixedbc} must be prescribed, namely.
\begin{equation}
\label{statbc}
\begin{array}{ll}
-\nabla^2 {\bf v} - \nabla [(\varepsilon_{min}^{-1} \varepsilon - 1) \nabla \cdot {\bf v}] = {\bf g} & \mbox{in } \Omega, \\
\left.
\begin{array}{l}
\partial_{n} {\bf e} \cdot {\bf n} + \nabla \cdot [\varepsilon_{min}^{-1} \varepsilon {\bf e}] - \nabla \cdot {\bf e} = 0, \\
\partial_{n} {\bf e}  \times {\bf n} = {\bf 0}.
\end{array}
\right\} & \mbox{on } \partial \Omega.
\end{array} 
\end{equation} 
This is actually the only real difference to be taken into account in order to extend to \eqref{eq2bis} the convergence analysis conducted in this paper. More precisely the final result that can be expected to hold for the fully discrete analog of \eqref{eq2bis}, defined as \eqref{eq7} under the same assumptions, is an $O(h^{\mu})$ error estimate, where $\mu \in (0,1]$ is such that the solution of \eqref{statbc} 
belongs to $[{\bf H}^{1+\mu}(\Omega)]^3$ for every ${\bf g} \in [L^2_0(\Omega)]^3$.\\

Another issue that is worth a comment is the practical calculation of the term $(\nabla \cdot \varepsilon {\bf e}_h^k, \nabla \cdot {\bf v})$ in \eqref{eq7}. Unless $\varepsilon$ is a simple function such as a polynomial, it is not possible to compute this term exactly. That is why we advocate the use of the 
well-known trapezoidal rule do carry out these computations. At the price of small adjustments in some  
terms involving norms of $\varepsilon$, the thus modified scheme is stable in the same sense as \eqref{stability}. 
Moreover the qualitative convergence result \eqref{convergence} remains unchanged, provided we require a little more regularity from $\varepsilon$. %except for the resulting mesh-independent constant ${\mathcal C}$. 
We skip details here for the sake of brevity.\\       

To conclude the authors should point out that the analysis conducted in this work can be adapted to the two-dimensional case without any significant modification. Clearly enough the same qualitative results apply to this case under similar assumptions.
\bigskip

\noindent \textbf{\underline{Acknowledgment:}} The second author gratefully acknowledges the financial support provided by CNPq/Brazil through grant 307996/2008-5, 
while part of this work was accomplished. \rule{2mm}{2mm} 

%\newpage

\end{document}